\documentclass[12pt]{amsart}

\usepackage{amsfonts}
\usepackage{amssymb,amsmath,amsthm}
\usepackage{verbatim}

\theoremstyle{definition}

\newtheorem{defi}{Definition}[section]
\theoremstyle{theorem}

\newtheorem{thm}{Theorem}[section]
\newtheorem{prop}{Proposition}[section]
\newtheorem{cor}{Corollary}[section]
\newtheorem{lm}{Lemma}[section]
\theoremstyle{remark}

\newtheorem{rem}{Remark}

\begin{document}
\title{Quaternionic contact structures in dimension 7 }
\author{David Duchemin}
\maketitle

\begin{abstract}

  The conformal infinity of a quaternionic-K\"ahler metric on a $4n$-manifold with boundary is a codimension $3$-distribution on the boundary called quaternionic contact. In dimensions $4n-1$ greater than $7$, a quaternionic contact structure is always the conformal infinity of a quaternionic-K\"ahler metric. On the contrary, in dimension $7$, we prove a criterion for quaternionic contact structures to be the conformal infinity of a quaternionic-K\"ahler metric. This allows us to find the quaternionic-contact structures on the $7$-sphere close to the conformal infinity of the quaternionic hyperbolic metric and which are the boundaries of complete quaternionic-K\"ahler metrics on the $8$-ball. Finally, we construct a $25$-parameter family of $Sp(1)$-invariant complete quaternionic-K\"ahler metrics on the $8$-ball together with the $25$-parameter family of their boundaries.

\end{abstract}

\section{Introduction}



In this paper we solve a boundary problem for quaternionic-K\"ahler metrics. This problem is a degenerate version of a problem initially posed for Einstein metrics. If $g$ is a metric on a manifold $M$ with boundary $N$, and $[b]$ is a conformal class of metrics on $N$, $[b]$ is the conformal infinity of $g$ if there exists a function $\rho$ positive in $M$ and vanishing to first order on $N$ such that $\rho^2 g $ extends continuously on $N$ with $\rho^2 g|_{T \mathbb{S}^3} \in [b]$. The standard example is the hyperbolic metric $g_{hyp}$ on the ball $B^{n+1}$ given by 
$$ g_{hyp} = 4 \frac{euc}{\rho^2} \, , $$
where $euc$ is the Euclidean metric on $\mathbb{R}^{n+1}$ and $\rho(x) = 1 - |x|^2$. The conformal infinity of $g_{hyp}$ is the conformal class of the round metric on $\mathbb{S}^n$.

The problem of finding complete Einstein metrics with prescribed conformal infinity on the ball was solved by Graham and Lee in \cite{Gra91}. In dimension $4$, one can search for selfdual Einstein metrics. LeBrun \cite{LeB82} shows using twistor theoretic arguments that a conformal metric on a $3$-manifold $N$ is the conformal infinity of a selfdual Einstein metric defined near $N$. However, a conformal metric on the sphere $\mathbb{S}^3$ is not always the conformal infinity of a complete selfdual Einstein metric on the ball $B^4$, see \cite{Biq02}. 

In the same way, the degenerate version is modeled on the quaternionic hyperbolic metric. Let $\mathbb{H}$ be the skew field of quaternions and $\mathbb{H}^n$ the $n$-dimensional $\mathbb{H}$-vector space. The action of the standard basis $(i,j,k)$ of imaginary quaternions gives endomorphisms $(I_1,I_2,I_3)$ of $\mathbb{H}^n \simeq \mathbb{R}^{4n}$. Each $I_i$ is an almost complex structure on $\mathbb{H}^n$ and one has the commutations rules $I_1 I_2 = - I_2 I_1 =I_3$. A such triple of endomorphisms on a real vector space $V$ is called a quaternionic structure on $V$. The quaternionic hyperbolic metric on the ball $B^{4n} \subset \mathbb{H}^n $ is given by$$ g_{\mathcal{H}} = \frac{euc}{\rho } + \frac{1}{ 4 \rho ^2 } ( (d \rho )^2 + ( I_1 d \rho )^2 + ( I_2  d \rho )^2 + ( I_3  d \rho )^2 ) \, , $$ where $\rho = 1 - |x|^2$ and $euc$ is the Euclidean metric. In this case, the function $\rho$ is positive in $B^{4n}$, vanishes to first order on $\mathbb{S}^{4n-1}$, and $[\rho^2 g_{\mathcal{H}}|_{T\mathbb{S}^{4n-1} } ]$ is a conformal class of degenerate metrics on $\mathbb{S}^{4n-1}$ with kernel 
$$ H^{can} = \cap_{i=1}^3 I_i d \rho |_{T \mathbb{S}^{4n-1 } } \, .$$
 The distribution $H^{can}$ is a so called quaternionic contact structure (\cite{Biq00} and \cite[p. 115]{Mont02}) whose definition in dimension $7$ is:

\begin{defi}
Let $H$ be an oriented distribution of codimension $3$ on a $7$-dimensional manifold $N$ and let $\mathcal{I}$ be the set of one forms vanishing on $H$. The distribution $H$ is called a quaternionic contact structure if $$\Lambda^2_+ H^*_x = \{ d \eta |_{H_x}\,  , \; \eta \in \mathcal{I} \} $$ is a rank three subbundle of $\Lambda^2 H^* $ such that the restriction to $\Lambda^2_+ H^*$ of the exterior product
 $$ \Lambda^2 H^* \otimes \Lambda^2 H^*  \rightarrow \Lambda^4 H^* \stackrel{\simeq}{\longrightarrow} \mathbb{R}  $$ gives a positive definite metric on $\Lambda^2_+ H^*$.
\end{defi}
If $H$ is a quaternionic contact structure in dimension $7$, a classical fact in $4$-dimensional linear algebra gives the existence of a unique conformal class $[g]$ of metrics on $H$ such that $\Lambda^2_+ H^*$ coincides with the space of selfdual $2$-forms with respect to $[g]$. Moreover, taking a local oriented orthonormal basis $(\frac{1}{\sqrt{2}} w_i  = \frac{1}{\sqrt{2}} d \eta_i|_{H}  )$ of $\Lambda^2_+ H^*$ with respect to a particular choice of metric $g$ in this conformal class, one gets a quaternionic structure $(I_i)_{i=1,2,3}$ on $H$ satisfying $w_i ( \cdot , \cdot ) = g( I_i \cdot , \cdot ) $ and defined up to a rotation by an element of $SO(3)$.

This description shows the link with the following definition given by Biquard in \cite{Biq00}: a quaternionic contact structure is a distribution $H$ of codimension $3$ on a manifold $N^{4n+3}$, locally given by three $1$-forms $(\eta_1, \eta_2, \eta_3)$ such that there exists a metric $g$ on $H$ and a quaternionic structure $(I_i)$ on $H$ satisfying the conditions $d \eta_i |_H = g ( I_i \cdot , \cdot )$. The conformal class $[g]$ is uniquely determined by $H$.

Our definition enlights the fact that in dimension $7$, quaternionic contact distributions form an open set in the set of codimension $3$ distributions. This fact is no more true in higher dimensions. 

Let us now come back to quaternionic-K\"ahler geometry. First, using the previous notations, we give the following definition:
\begin{defi}
A metric $g$ on a manifold $M$ with boundary $N$ is asymptotically quaternionic hyperbolic (AQH) if one has a quaternionic contact structure $H$ on $N$ with compatible metric $g_H$ on $H$ and a function $\rho$, positive in $M$ vanishing to first order on $N$ such that on a neighbourhood $] 0 , a ] \times N$ of $N$, the behaviour of $g$ near $N$ is given by 
$$g \sim  \frac{1}{\rho^2} ( d \rho ^2 +  \eta_1 ^2 + \eta_2^2 + \eta_3^2)  + \frac{1}{\rho} g_H \; \mbox{when}  \; \rho \rightarrow 0 . $$
 The quaternionic contact structure $H$ is called the conformal infinity of $g$. If $g$ is also quaternionic-K\"ahler, one says that $g$ is asymptotically hyperbolic quaternionic-K\"ahler (AHQK).
\end{defi}

 Biquard \cite{Biq00} has shown that any quaternionic contact structure of dimension $4n+3 \geq 11$ is at least locally the conformal infinity of a unique AHQK metric. Moreover, he showed in \cite{Biq02} that a quaternionic contact structure on $\mathbb{S}^{4n+3}$ with $4n+3 \geq 11$ and close to the canonical one is the conformal infinity of a AHQK complete metric on the ball $B^{4n+4}$. The question remains open in dimension $7$. 

In this paper, we answer this last question. We show that the conformal infinity of an AHQK $8$-manifold must satisfy an additionnal integrability property which is empty in higher dimensions. Conversely, we prove that an integrable quaternionic contact $7$-manifold is the conformal infinity of a unique AHQK manifold.
\begin{defi}
Let $H$ be a quaternionic contact structure on a manifold $N$ of dimension $7$ and choose a compatible metric $g$. The quaternionic contact structure $H$ is called integrable if there exists a local oriented orthonormal basis $(d \eta_i|_H )$ of $\Lambda^2_+ H^* $ and vector fields $(R_1, R_2,R_3)$ satisfying
\begin{itemize}
\item  $i_{R_i}\eta_j = \delta_{ij}$,
\item $i_{R_i} d \eta_j |_H = - i_{R_j} d \eta_i |_H $.
\end{itemize}
This property does not depend on the choice of metric $g$ inside the conformal class.
\end{defi}
We can now give the statements of the main results. 
\begin{thm} \label{thm1}
Let $H$ be a real analytic quaternionic contact structure on a manifold $N^7$. Then $H$ is the conformal infinity of an AHQK metric $g$ defined on a neighbourhood of $N$ and admitting a real analytic extension on the boundary with pole of order 2 iff $H$ is integrable. Moreover, the germ of $g$ along $N$ is uniquely determined by $H$.
\end{thm}

Using \cite{Biq02} and this theorem, we can fill in the $8$-ball by globally defined complete AHQK metrics whose boundaries are close to the canonical quaternionic contact structure $H^{can}$:

\begin{cor} \label{cor}
Let $H$ be an integrable quaternionic contact structure on $\mathbb{S}^7$, close to the canonical distribution $H^{can}$. Then $H$ is the conformal infinity of a complete AHQK metric on the ball $B^8$.
\end{cor}

Among the integrable quaternionic contact structures on $\mathbb{S}^7$, we show the existence of an interesting family of $Sp(1)$-invariant integrable quaternionic contact structures on the $7$-sphere:
 
\begin{thm} \label{thm2}
Let $H^{can}$ be the canonical quaternionic contact structure of $\mathbb{S}^7$. Let $\mathcal{H}$ be the set of integrable $Sp(1)$-invariant quaternionic contact structures and $\mathcal{G}$ be the group of diffeomorphisms of $\mathbb{S}^7$ commuting with the $Sp(1)$-action. There is a neighbourhood $\mathcal{V}$ of $[H^{can}]$ in $\mathcal{H}/ \mathcal{G}$ which is homeomorphic to the quotient of a $35$-dimensional ball $B^{35}$ by the isotropy group $Sp(2)$ of $H^{can}$. One obtains a $25$-parameter family of integrable quaternionic contact structures.
\end{thm}

Then, we can construct a family of $Sp(1)$-invariant complete quaternionic-K\"ahler metrics on the $8$-ball :

\begin{cor} \label{cor2}
Let $g_{\mathcal{H}}$ be the quaternionic hyperbolic metric on the $8$-ball. There exists a $25$-parameter family of $Sp(1)$-invariant AHQK metrics with boundaries close to the boundary $H^{can}$ of $g_{\mathcal{H}}$. 
\end{cor}

This examples generalize a $3$-parameter family constructed by Galicki in \cite{Gal91}. These metrics are obtained by quaternionic quotient of the hyperbolic quaternionic space $\mathbb{H} \mathcal{H}^3$ and all have isometry group strictly greater than $Sp(1)$.

The paper is organized as follows. In section 2, we
construct a connection associated to each compatible metric. A part $T^W$ 
of its torsion gives a conformal invariant named vertical torsion. The vanishing of $T^W$ is equivalent to the integrability of $H$.

In the third section, we study the boundaries of AHQK manifolds and we show that they are integrable. This gives the motivation to study more carefully the torsion and the
curvature of this case. In particular, the curvature on $H$ looks like that of anti-selfdual Riemannian $4$-manifolds except for an additional term coming from the Bianchi identity. The computation is done in section
4.

 Still assuming the integrability condition, we construct an integrable CR-manifold, the twistor space of the quaternionic contact structure. This is done in section 5 and gives the converse statement to the third section, namely that a quaternionic contact structure with vanishing vertical torsion is the boundary of a unique AQH manifold of dimension $8$. 

Section $6$ is devoted to the study of deformations of $H^{can}$. Then, we describe in detail the case of $Sp(1)$-invariant deformations of the $7$-sphere and show the existence of a $25$-parameter family of integrable $Sp(1)$-invariant deformations of $H^{can}$.

\bigskip

\textbf{Acknowledgments}: This paper is a part of the author's doctoral thesis; in this connection thanks are due to O. Biquard for his extremely helpful comments.

\section{Construction of the connection}

In the following, one has a smooth manifold $N$ of dimension $7$, a quaternionic contact structure $H$ on $N$ and $g$ a fixed compatible metric $g$ on $H$. We fixe local contact forms $(\eta_1, \eta_2, \eta_3)$ and a local quaternionic structure $(I_i)$ on $H$ such that $d\eta_i ( \cdot , \cdot ) = g ( I_i \cdot , \cdot )$ on $H$.

In the first three parts of this section, we construct an adapted connection associated to $g$. This connection will be used in the twistorial construction of section 5. To look at the conformal invariance of this twistorial construction, we will need to know how a conformal change of metric changes the connection. This is done in part 5 of this section.

\subsection{Partial connection}
If $N$ is a manifold, $E$ a vector bundle and $D$ a distribution on $N$, a $D$-connection on $E$ is a differential operator 
\begin{equation*}
\nabla : \Gamma(E) \rightarrow \Gamma ( D^* \otimes E ) \, , 
\end{equation*}
satisfying the Leibniz rule $\nabla ( f \, s ) = (df)|_D \otimes s + f
\nabla s $ for every function $f$ and section $s$ of $E$.

\begin{lm}
\label{partial} Assume that $W$ is a distribution on $N$ giving a splitting $TN = H \oplus W$. There exists a unique $H$-connection $\nabla $ on $H$
preserving the metric $g$ and such that the torsion satisfies 
\begin{equation*}
\forall X, Y \in H, \; ( T_{X, Y} )_H = 0 \, , 
\end{equation*}
where the subscript $H$ indicates the projection on $H$ in the direction of $W$.\end{lm}

\begin{proof}
 If $\nabla$ is such a connection, we must have for every sections $X$, $Y$ and $Z$ of $H$ the Koszul formula 
$$\begin{array}{rcl}
2 g ( \nabla_X Y , Z ) & = & X. g( Y, Z )+ Y. g ( Z, X) - Z. g (X,Y ) \\  
& & + g ( [ X, Y ] _H, Z ) - g ( [X, Z]_H, Y ) - g( [ Y, Z ] _H , X ) \, .
\end{array} $$
It gives both uniqueness and existence.
\end{proof}

 Otherwise stated, the vector fields $X$, $Y$, $Z$ are sections of $H$, and given a complementary $W$, a vector field $R$ is a section of $W$, and $(R_1,R_2,R_3)$ is the dual basis of $%
(\eta_1|_W, \eta_2|_W, \eta_3|_W)$. We equip $W$ with the metric $\sum_i \eta_i^2 $.
\begin{rem}
If $W$ is a complement to $H$, the torsion of the $H$-connection
associated to $W$ on $H$ satisfies
\begin{equation*}
T_{X,Y} = -[X,Y]_W = \sum_{i=1}^3 d\eta_i ( X, Y ) R_i \, .
\end{equation*}
\end{rem}

\subsection{Extension of the connection}

\begin{lm} \label{1}
Let $W$ be a complement of $H$ in $TN$. One can find a unique connection $%
\nabla^W$ on $N$ such that :

\begin{itemize}
\item[(i)]  $\nabla^W$ preserves the splitting $TN= H \oplus W$ and the
metrics on $H$ and $W$,

\item[(ii)] if $X,\, Y \in H$ and $R,\, R' \in W$, then $(T_{X,Y})_H =0$ and $(T_{R,R'} )_W = 0$,

\item[(iii)]  the torsion $T$ satisfies 
\begin{eqnarray}
& \forall X \in H \; , \; T_X^W := ( R \mapsto (T_{X, R})_W)  \in \mathfrak{so}%
(W)^{\perp} \, , & \\
& \forall R \in W \; , \; T_R ^H := ( X \mapsto (T_{R,X})_H) \in \mathfrak{so}(H)^{\perp}
\, , &
\end{eqnarray}
\end{itemize}
\end{lm}

\begin{proof}
Let $\nabla$ be the partial connection on $H$ defined by lemma \ref{partial}. We extend it to a true connection which preserves the metric on $H$, still denoted by $\nabla$. If $a \in \Gamma ( W^* \otimes \mathfrak{so}(H))$, the connection $\nabla' = \nabla + a$ is metric and its torsion $T'$ satisfies 
$$T'_{R,X}  = \nabla'_R X - [ R, X ] _H = T_{R,X} + a_R ( X) \, , $$
so that there exists a unique $a_R$ which annihilates the $\mathfrak{so}(H)$-part of $T_{R,\cdot}$. The covariant derivatives in the direction of $W$ are defined in the same way.
\end{proof}

We put $\alpha_{ij} ( X) = d \eta_j ( R_i , X )$. One has 
\begin{equation*}
T^W_X ( R_i) = \nabla_X^W ( R_i ) - [ X, R_i ]_W = \nabla_X ^W R_i -
\sum_{j=1}^3 \alpha_{ij} (X ) R_j \, , 
\end{equation*}
from which we obtain
\begin{equation*} 
\nabla_X^W R_i = -\frac 12 \sum_{j=1}^3 ( \alpha_{ji} (X) - \alpha_{ij} (X) )
R_j
\end{equation*}
and
\begin{equation} \label{torsion}
 T_X^W (R_i ) =- \frac 12 \sum_{j=1} ^3 ( \alpha_{ji} (X)
+ \alpha_{ij} ( X) ) R_j \, . 
\end{equation}

\subsection{Reducing torsion}
\label{deux}  We search now a particular choice of $W$ giving the simplest torsion. To fix the
notations, we recall some basic facts about representations of $SO(4)$.

The universal covering of $SO(4)$ is $Spin(4) = Sp(1)\times Sp(1) $ where $
Sp(1)$ is the group of unitary quaternions. Let $S_+$ and $S_-$ be the
representations of the first and the second factor respectively on $\mathbb{H%
} \simeq \mathbb{C}^2$. The irreducible representations of $Spin(4)$ are the 
$S_+^m \otimes S_-^n$ where $S_+^m$ and $S^n_-$ are the symmetric power of
order $m$ and $n$ of $S_+$ and $S_-$ respectively. The following
Clebsch-Gordan formula gives the irreducible decomposition of tensorial
products~: 
\begin{equation*}
S_+^n \otimes S_+^p \simeq S_+^{n+p} \oplus S_+^{n+p-2} \oplus \cdots \oplus
S_+^{n-p} \, , \; p \leq n \, . 
\end{equation*}

The real irreducible representations of $SO(4)$ are the real parts of $%
S_+^n \otimes S_-^m$ with $n+m$ even. We will denote them by $S^{n,m}$. In
particular, we have 
\begin{equation*}
\mathbb{R}^4 \simeq S^{1,1} \, , \; \Lambda ^2_+ \simeq S^{2,0} \, , \;
\Lambda^2_- \simeq S^{0,2} \, . 
\end{equation*}

We now give the explicit isomorphism $\mathbb{R}^4 \otimes Sym^2 (
\Lambda ^2_+ ) \simeq S^{5,1} \oplus S^{3,1} \oplus S^{1,1} $. Let $(I_1, I_2,
I_3 )$ be a quaternionic structure on $\mathbb{R}^4$ given a $SO(3)$-trivialization of $\Lambda^2_+$. Then 
\begin{equation*}
\begin{array}{rcl}
S ^{5,1} & \simeq & \{ \sum_{i,j} a_{ij} \otimes I_i \otimes I_j \, , \;
a_{ij} =a_{ji }\in \mathbb{R}^4 \; \; and \; \; \forall j \, , \; \sum_i I_i
a_{ij} = 0 \} \, , \\ 
S^{1,1} & \simeq & \{ \sum_i r \otimes I_i \otimes I_i \, , \; r \in \mathbb{%
R}^4 \} \, , \\ 
S^{3,1} & \simeq & \{ \sum_{i,j} ( I_i r_j + I_j r_i ) \otimes I_i \otimes
I_j \,, \; r_i \in \mathbb{R}^4 \; , \; \sum_i I_i r_i = 0 \} \, .
\end{array}
\end{equation*}

In our case we have the natural identification 
$$ W \simeq \Lambda^2_+ H^* \, , \; R_i \mapsto d \eta_i|_H$$ so that $T^W$ becomes a section of $H^* \otimes End
( \Lambda^2_+ H ^* ) $. We put $w_i = d\eta_i |_H$ and $w_i^*$ the dual basis.

\begin{rem}
The metric $g$ allows us to identify $H^*$ and $H$ and we use it throughout
the text. In particular $\Lambda^2_+ H^* $ can be considered as a subspace of
the space of 2-forms or that of skew-symmetric endomorphisms.
\end{rem}

\begin{prop}
For each choice of compatible metric $g$ on $H$, there is a
unique complement $W^g$ of $H$ such that $T^{W^g} \in \Gamma( S^{5,1} )$.
\end{prop}

\begin{proof} Let $W$ be transverse to $H$ and $(R_1, R_2, R_3 )$ be the dual basis of $(\eta_1, \eta_2, \eta_3 )$ on $W$. We have obtained in (\ref{torsion})
$$T^W = - \frac 12 \sum_{j=1}^3 ( \alpha_{ij} + \alpha_{ji} ) \otimes w_i^{*} \otimes w_j \, , $$

If $W'$ is another complementary to $H$ spanned by the vectors $R_i' = R_i + r_i$ with $r_i \in H$, then $\alpha_{ij}' = i_{R_i'} d \eta_j|_H = \alpha_{ij} + ( I_j r_i)^{\flat} $ ($\flat$ and $\sharp$ are the usual musical isomorphisms). With the explicit decomposition of $H^* \otimes Sym^2 ( \Lambda^2_+ H) $ we wrote down, the existence and the uniqueness of $W$ follow easily.
\end{proof}

\begin{rem} Another choice of complementary does not change the $S^{5,1}$ part of the torsion.

\end{rem}

\subsection{Derivation of the quaternionic structure}
\label{horder}
We fix $W = W^g$ and note $\nabla$ the corresponding connection. This
connection is metric and so preserves the bundle $\Lambda^2_+ H^*$ : 
\begin{equation*}
\nabla I_j|_H = \sum_{i=1}^3 \gamma_{ij}\otimes I_i \, . 
\end{equation*}
Here we just look at the derivation in the direction of $H$, i.e. $\gamma_{ij}
\in H^*$.

Let $X, Y, Z \in \Gamma (H )$, and $\mathfrak{a}$ be skew-symmetrisation in $X$
$Y$ and $Z$. One has the identity 
\begin{equation} \label{detai}
d ( d \eta_j ) (X,Y, Z ) = \mathfrak{a}( \nabla d \eta_j ) ( X, Y ,Z ) + d
\eta_j ( T_{X,Y} , Z ) + d \eta_j ( T_{Z, X} , Y ) + d \eta_j ( T_{Y, Z} , X
) \, , 
\end{equation}
which can be rewritten in the following form : 
\begin{equation*}
\sum_{i}( g ( \gamma_{ij} ^{\sharp} + \alpha_{ij} ^{\sharp }, X ) I_i + I_i
X \wedge ( \gamma_{ij}^{\sharp} + \alpha_{ij} ^{\sharp } ) ) = 0 \, . 
\end{equation*}
Projecting on $\Lambda^2_+ H$ and $\Lambda^2_- H$ gives the equivalent
condition 
\begin{equation*}
\forall j \in \{1,2,3 \} \, , \; \sum_{i=1} ^3 ( \alpha_{ij } +\gamma_{ij} )
\circ I_i = 0 \, .
\end{equation*}
But our particular choice of complementary vector bundle ensures that $%
\sum_i ( \alpha_{ij} + \alpha_{ji} ) \circ I_i = 0$, hence we get 
\begin{equation*}
\gamma_{ij} = - \frac 12 ( \alpha_{ij} - \alpha_{ji} ) \, . 
\end{equation*}

\subsection{Conformal change}

 Let $\eta^{\prime}= f^2 \eta $ be such a conformal change, and $(R_1,R_2,R_3)$ be the dual basis of $(\eta_1, \eta_2, \eta_3)$ on $W^g$. We put $(R^{\prime}_1, R_2^{\prime},R_3^{\prime})$ the dual basis of $(f^2 \eta_1 ,f^2 \eta_2 , f^2 \eta_3 ) $ on $W^{f^2g}$.

\begin{prop}
\label{conforme} The conformal change of metric corresponds to the following change of basis of associated complementaries : 
\begin{equation*}
R_i^{\prime}= f^{-2} (R_i +r_i) \, , 
\end{equation*}
where $r_i^{\flat} = 2 f^{-1} df|_H \circ I_i$ (the musical isomorphisms $%
\sharp$ and $\flat$ are taken with respect to $g$ on $H$, after restriction
if necessary for 1-forms). Moreover, we get 
\begin{equation*}
i_{R^{\prime}_i} d \eta_j^{\prime}|_H + i_{R^{\prime}_j} d
\eta_i^{\prime}|_H = i_{R_i} d \eta_j |_H +i_{R_j} d \eta_i |_H \, . 
\end{equation*}
\end{prop}

\begin{proof}
We put $\alpha_{ij}' = i_{R_i'} d\eta_j' |_H$. We have $\eta'_i(R_j') = f^2 \eta_i ( R_j') = \delta_{ij}$ so that $R_i'= f^{-2}( R_i + r_i)$ with $r_i \in H$ and finally
\begin{eqnarray}
 & \alpha_{ij}' + \alpha_{ji}' = \alpha_{ij} + \alpha_{ji} -(r_i^{\flat} \circ I_j +r_j^{\flat} \circ I_i ) - 4  \delta_{ij} f^{-1} df|_H \, . &
\end{eqnarray}
 
The conformal change left $S^{5,1}$, $S^{3,1}$ and $S^{1,1}$ globally invariant and $(r_i^{\flat} \circ I_j +r_j^{\flat} \circ I_i ) + 4 \delta_{ij} f^{-1} df|_H \in S^{3,1} \oplus S^{1,1} $ therefore the conditions $\alpha_{ij}' + \alpha_{ji}' \in S^{5,1} $ and $\alpha_{ij} + \alpha_{ji} \in S^{5,1}$ imply $(r_i^{\flat} \circ I_j +r_j^{\flat} \circ I_i ) + 4 \delta_{ij} f^{-1} df|_H =0 $ and the lemma follows. 
\end{proof}

\begin{cor}
The torsion $T^{W^g}$ associated to the Carnot-Carath\'eodory metric is
conformally invariant. We call it the vertical torsion and denote it by $%
T^{W^g}$ or $T^W$.
\end{cor}

\begin{proof}
If we change the metric in the conformal class, the $2$-forms $w_i$ are multiplied by the conformal factor and elements of the dual basis are multiplied by its inverse. So the only thing we must look at is the invariance of $( \alpha_{ij} + \alpha_{ji})_{i,j} $ which follows from \ref{conforme}.

\end{proof}

Let us summarize the results we have obtained in the following proposition.

\begin{prop} \label{ifint}
 Let $(N,H)$ be a quaternionic contact structure. The integrability of $H$ does not depend on the choice of an adapted metric on $H$. Moreover, if $g$ is particular a choice of compatible metric on $H$, the following conditions are equivalent :
\begin{itemize}
\item The distribution $H$ is integrable.
\item The torsion $T^{W^g}$ vanishes.
\item For any choice of complementary distribution $W$, the $S^{5,1}$ part of the torsion vanishes.
\item For any choice of oriented orthonormal basis $(\frac{1}{\sqrt{2}}d\eta_i |_H)$ of $H^+$ and any choice of vector fields $(R_1,R_2,R_3)$ such that $$\eta_j(R_i) =\delta_{ij} \, , $$
the $S^{5,1}$ part of $(i_{R_i} d \eta_j |_H + i_{R_j} d \eta_i |_H )_{i,j}$ vanishes.
\end{itemize}

\end{prop}

In the study of the twistor space, we will need to know how the connection
is changed when the metric is multiplied by a conformal factor. We put $%
\theta = f^{-1} d f $. Recall that we write $\theta^{\sharp} $ for $(
\theta|_H)^{\sharp}$ and that the change of complementary distribution is
parametrized by $R_i' = f^{-2} ( R_i -2 I_i \theta^{\sharp} ) $. The following lemma will be useful in the twistorial construction.

\begin{lm}
\label{changeconnection} Following the notations of 2.2, the new connection $%
\nabla^{\prime}$ is given by 
\begin{equation*}
\begin{array}{rcl}
\nabla_X ^{\prime} & = & \nabla_X + \theta(X) + \theta^{\sharp} \wedge X +
\sum_i I_i \theta^{\sharp} \wedge I_i X+ \sum_i\langle I_i \theta^{\sharp}, X \rangle
I_i  \\ 
\nabla_{R_i} ^{\prime} & = & \nabla_{R_i} + \theta(R_i) + 2 |
\theta^{\sharp} | ^2 I_i+ 2 \theta^{\sharp} \wedge I_i \theta^{\sharp} -
\frac 12 \sum_j (\alpha_{ij}^{\sharp} +\alpha_{ji} ^{\sharp} ) \wedge I_j
\theta^{\sharp} +2 ( I_i \nabla \theta^{\sharp} )^{\mathfrak{so}(H)}
\end{array}
\end{equation*}
\end{lm}

\begin{proof}
 We put $\nabla ' = \nabla +\theta + a $ and $\nabla^1 = \nabla + \theta$. The connection $\nabla^1$ preserves $f^2 g$ and its torsion is
$$\begin{array}{rcl}
T^1_{X,Y} & = & \sum_i d \eta_i ( X, Y ) R_i + \theta(X) Y - \theta(Y) X \\
 & = & T'_{X, Y} - \sum_i d\eta_i ( X, Y, ) r_i + \theta(X) Y - \theta(Y) X \\
\end{array} $$
so that $a_X Y - a_Y X = \sum_i d \eta_i ( X, Y) r_i - \theta(X) Y + \theta (Y) X $. The connections $\nabla'$ and $\nabla^1$ both preserve $f^2 g $ hence $a$ is a 1-form with values in $\mathfrak{so} (H)$. The skew-symmetrisation in the two first variables gives an isomorphism $ H^* \otimes \mathfrak{so}(H) \rightarrow \Lambda ^2 H^* \otimes H$, with inverse $b$
$$ \langle b ( c)_X Y , Z \rangle = \frac 12 \left ( \langle c(X,Y) ,Z \rangle + \langle c (Z,X) ,Y \rangle- \langle c(Y, Z ) , X \rangle \right ) \,  $$
 from which we deduce the first part of the lemma.

We now look at the change of the connection in the direction of $W^g$. If $U \in TS$, $U_{V/W}$ is its projection on $V$ in the direction of $W= W^g$. We have 
$$\begin{array}{rcl}
a_{R_i} X & = & \nabla_{R_i} ' X  - \nabla_{R_i } X- \theta(R_i) X \\
 & =& \nabla_{R_i + r_i} ' X - \nabla_{R_i} X - \theta(R_i) X - \nabla_{r_i} ' X \\ 
\end{array} $$
 Introducing the torsion, we obtain
$$\begin{array}{rcl}
a_{R_i} X &= & (T_{R_i + r_i, X  })_{V/W'} - (T_{R_i, X} )_{V/W} +[R_i, X]_{V/W'} - [R_i, X]_{V/W} \\
 & &  - \theta(R_i) X - \nabla_{r_i} ' X +[r_i, X ] _{V/W'}   \\
 & = &  (T_{R_i + r_i, X  })_{V/W'} - (T_{R_i, X} )_{V/W} + \sum_j d \eta_j ( R_i, X ) r_j  \\
 & & - \theta(R_i) X - \nabla_{X}' r_i
\end{array} $$
But $a_{R_i} \in \mathfrak{so}(H)$, so that it suffices to compute the skew-symmetric part of the right term in the previous equality. The contributions of the torsions vanish by definition, that of $\sum_j d\eta_j (R_i, X) r_j $ is 
$$\frac 12 \sum_j \alpha_{ij}^{\sharp} \wedge r_j = - \sum_j \alpha_{ij}^{\sharp} \wedge I_j \theta^{\sharp}\, , $$
and that of $\nabla' r_i $ is
$$- 2 \theta^{\sharp} \wedge I_i \theta^{\sharp}  - 2 | \theta^{\sharp} |^2 I_i  + ( \nabla r_i ) ^{\mathfrak{so}(H)}$$
where the exponent $\mathfrak{so}(H)$ means the orthogonal projection on $\mathfrak{so}(H)$. But using part 2.4, we get
$$ \begin{array}{rcl}
\nabla r_i & = & - 2 \nabla ( I_i \theta^{\sharp})  = - 2 (\nabla I_i ) \theta^{\sharp} - 2 I_i \nabla \theta^{\sharp} \\
 & =& \sum_j ( \alpha_{ji} - \alpha_{ij} ) \otimes I_j \theta^{\sharp} - 2 I_i \nabla \theta^{\sharp} \, . 

\end{array} $$
Mixing all this together gives the lemma.
\end{proof}

\subsection{Higher dimensional case} 

Let us do some remarks about what is going on in higher dimensions. Let $H$ be a quaternionic contact structure on a manifold $N^{4n+3}$ with $n>1$ and $g$ be a compatible metric on $H$. In the same way and always with the same notations, one can show that there exists a unique complementary $W^g$ such that  
$$\sum_i ( \alpha_{ij} + \alpha_{ji} ) \circ I_i =0 $$
for all $j$. 

On the other hand, \ref{partial} is always true and give a metric $H$-connection $\nabla$ on $H$. Then, using (\ref{detai}) one can show that in fact $\alpha_{ij} + \alpha_{ji} = 0$ and that $\nabla$ preserves not only the metric but also the $Sp(n)Sp(1)$ structure on $H$. Hence, there is no integrability condition. It is the reason why all quaternionic contact structures in dimension strictly greater than $7$ are the boundaries of AHQK metrics.

\section{Conformal infinity of AHQK manifolds}

In this section, we will study the conformal infinity of an AQH quaternion-K\"ahler manifold. We find a particular trivialization
of the quaternionic structure admitting an analytic extension to the
boundary with pole of order 2. Then, we use it to show that the quaternionic contact structure on the boundary is integrable.

\subsection{Twistor space and asymptotic development} \label{3.1}

The following is essentially the work of \cite[III.2]{Biq00} and \cite{LeB91}%
. Let $(M,g)$ be an AHQK manifold
of dimension $8$ and suppose that the metric $g$ admits an analytic
extension to the boundary $N$. We will apply the twistor machinery to obtain
a particular choice of local trivialization of the quaternionic structure in
a neighbourhood of the boundary. 
The twistor space \cite{Sal82} of $M$ is a $5$-dimensional holomorphic
manifold with the following data :

\begin{itemize}
\item  a holomorphic contact structure $\eta$ with values in a line bundle $L
$;
\item  a family of dimension $8$ of compact genus zero curves $(\mathcal{C}%
_m)_{m \in M^{\mathbb{C}}}$ with normal bundle $\mathcal{O}(1) \otimes 
\mathbb{C}^4$;
\item  an hypersurface $N^{\mathbb{C}} \subset M^{\mathbb{C}}$ of curves tangent to the contact distribution;
\item  a compatible real structure $\tau$, without fixed points.
\end{itemize}
\begin{rem}
 $M$ is the real slice of $M^{\mathbb{C}}-N^{\mathbb{C}}$ and $N$ that of $N^{\mathbb{C}}$.
\end{rem}

On each $\mathcal{C}_m$, the line bundle $L$ is isomorphic to $\mathcal{O}(2)$
so that $\mathcal{L}_m = H^0 ( \mathcal{C}_m , Hom( T \mathcal{C}_m , L ) )$
is a line bundle on $M^{\mathbb{C}}$. By restriction, the $1$-form $\eta$ gives a section 
$\Theta$ of $\mathcal{L}$ and $S^{\mathbb{C}}$ is the null set of $\Theta$.
We choose local square root $L^{1/2}$ of $L$, but the conclusions do not
depend on this choice. Let us define 
\begin{equation*}
\begin{array}{rcl}
E_m & = & H^0 ( \mathcal{C}_m, L^{-1/2} \otimes N_m )\, , \\ 
H_m & = & H^0 ( \mathcal{C}_m , L^{1/2} ) \, ,
\end{array}
\end{equation*}
so that
$$T_m M^{\mathbb{C}} = E_m \otimes H_m \, . $$
For $m \notin N^{\mathbb{C}}$ and $u, v \in H_m$, the Wronskian $w(u \wedge v ) = u dv - v du $ defines a two form 
\begin{equation}
w_H : \Lambda ^2 H^0 ( \mathcal{C}_m , L^{1/2} ) \overset{w}{\rightarrow} 
\mathcal{L}_m \overset{\Theta^{-1}}{\rightarrow} \mathbb{C} \, , 
\end{equation}
and therefore a $SO_3 ( \mathbb{C} )$-structure $w_H \otimes w_H$ on $H^0 ( 
\mathcal{C}_m , L ) \simeq Sym^2 ( H_m ) $.

The normal bundle $N_m$ of a curve $\mathcal{C}_m$ has a natural
identification with $\ker \eta$ if $m \notin N^{\mathbb{C}} $ so that we
have a well defined $2$-form 
\begin{equation*}
 \Lambda ^2 H^0 ( \mathcal{C}_m , N_m ) \overset{d \eta }{\rightarrow} Sym^ 2
(H_m ) \, . 
\end{equation*}
The choice of a $SO_3 ( \mathbb{C} )$-trivialization on $Sym^2 ( H_m ) $
exhibits three $2$-forms $w_1$, $w_2$, $w_3$ giving the $Sp_2(\mathbb{C})Sp_1(\mathbb{C})$ structure. The complexified quaternionic-K\"ahler metric is
\begin{equation} \label{metric}
 g = w_E \otimes w_H \; \mbox{on}\;  E_m \otimes H_m \, 
\end{equation}
where 
$$w_E : E_m \stackrel{d \eta}{\rightarrow} \mathbb{C} \, . $$
We now look at the contact structure on the boundary. Let $l : \mathcal{L}
\rightarrow \mathbb{C}$ be a local choice of trivialization of $
\mathcal{L}$ in a neighbourhood of $s \in N^{\mathbb{C}}$ and extend it on $%
M^{\mathbb{C}}$. In the same way, we obtain a symplectic form 
\begin{equation*}
\hat{w}_H : \Lambda^2 H^0 ( \mathcal{C}_m , L^{1/2} ) \overset{w}{\rightarrow} \mathcal{L}%
_m \overset{l}{\rightarrow} \mathbb{C} \, , 
\end{equation*}
and thus a $SO_3 (\mathbb{C})$-metric $\hat{w}_H \otimes \hat{w}_H = l^2
\Theta^2 w_H \otimes w_H$. We choose a local $SO_3 ( \mathbb{C})$%
-trivialization $Sym^2 ( H_m ) \rightarrow \mathbb{C}^3 \, . $

If $s \in N^{\mathbb{C}}$, one has $T \mathcal{C}_s \subset \ker \eta $ hence $\eta$
gives three 1-forms $(\eta_1, \eta_2, \eta_3 ) $ along $N^{\mathbb{C}}$ 
\begin{equation*}
H^0 ( \mathcal{C}_s, N_s ) \overset{\eta}{\rightarrow} H^0 ( \mathcal{C}_s,
L ) \simeq Sym^2 ( H_s ) \rightarrow \mathbb{C}^3 \, . 
\end{equation*}
On the other hand, on $M^{\mathbb{C}}- N^{\mathbb{C}}$ we obtain three
2-forms
\begin{equation*}
\wedge ^2 H^0 ( \mathcal{C}_m , N_m) \overset{d \eta}{\rightarrow } Sym^2 (
H_m ) \rightarrow \mathbb{C}^3 \, , 
\end{equation*}
which can be written as $l^2 \Theta^2 w_i$ with $w_i$ defining the
quaternionic structure of $M^{\mathbb{C}} - N^{\mathbb{C}} $.

We put $\rho = l \Theta : M^{\mathbb{C} } \rightarrow \mathbb{C} $.

\begin{lm}
The forms $w_i$ have pole of order $2$ along $N^{\mathbb{C}}$. More precisely, the
2-forms $l^2 \Theta^2 w_i$ are defined on $N^{\mathbb{C}}$ and satisfy  
\begin{equation*}
l^2 \Theta^2 w_i = - d\rho \wedge \eta_i + \frac 12 \sum_{r,s}
\varepsilon^{rsi} \eta_r \wedge \eta_s \, 
\end{equation*}
on $N^{\mathbb{C}}$ where $\varepsilon^{rsi}$ is the signature of the permutation $(r,s,i)$ of $(1,2,3)$.
\end{lm}

\begin{proof}
Because the $w_i$ define a quaternionic structure, we need only show that $i_{\partial/\partial \rho }l^2 \Theta^2 w_i = - \eta_i $ to obtain the lemma. We take $s \in N^{\mathbb{C}}$. 

There exists a section $\phi$ of $N_s$ along $\mathcal{C}_s$ such that $\eta(\phi) = 0$ and $i_{\phi} d\eta|_{T\mathcal{C}_s} \neq 0$, cf \cite[lemma III.2.5]{Biq00}. We normalize $\phi$ in order to have $l( i_{\phi} d\eta|_{T\mathcal{C}_s}) =1$. It is a vector in $T_s M^{\mathbb{C}}$ with the properties $d\rho ( \phi ) =1$ and $\eta_i ( \phi ) =0$. 
Remark that whereas the symplectic form $d \eta$ is not defined along $N^{\mathbb{C}}$, the 3-form $\eta \wedge d \eta $ admits an extension to $N^{\mathbb{C}}$. By restriction, we have
$$ \Theta d \eta = \eta \wedge d \eta \in H^0 ( \mathcal{C}_m , T^* \mathcal{C}_m \otimes \Lambda^2 N_m^* \otimes L^2 ) = \mathcal{L}_m \otimes \Lambda^2 T^* M^{\mathbb{C}} \otimes H^0 ( \mathcal{C}_m, L ) \, . $$

If $u$ is tangent to $\mathcal{C}_s$ and $\sigma \in H^0 ( \mathcal{C}_s , N_s )$, then
$$ \eta \wedge d \eta ( u, \phi , \sigma ) = \eta ( \sigma ) d \eta ( u, \phi) \, , $$
i.e. $ i_{\phi } l\Theta d \eta = - \eta $ and finally
$$ i_{\phi} l^2 \Theta^2 w_i  = -\eta_i \, . $$
\end{proof} 
The intersection of the kernels of $\rho^2 w_1$, $\rho^2 w_2$, and $\rho^2
w_3$ on $N^{\mathbb{C}}$ is $$H^{\mathbb{C}} = H^0 ( \mathcal{C}_s , N_s \cap \ker \eta \cap T 
\mathcal{C}_s^{\perp d \eta })$$ and coincides with the contact structure of
the boundary. The symplectic form $w_i$ has well defined terms of order $-1$ on $H^{%
\mathbb{C}}$ and one can show \cite[Lemma III.2.6]{Biq00} that 
\begin{equation*}
w_i = \frac{1}{\rho^2} w_{i,-2} + \frac{1}{\rho} w_{i, -1} + \cdots \, , 
\end{equation*}
with 
\begin{equation*}
w_{i,-2} = - d\rho \wedge \eta_i + \frac 12 \sum_{r,s} \varepsilon^{rsi}
\eta_r \wedge \eta_s \, , \; w_{i, -1} |_{H^{\mathbb{C}}} = d\eta_i |_{H^{%
\mathbb{C}}} \, . 
\end{equation*}
If we put 
$$\hat{E}_m = H^0(\mathcal{C}_s, (N_s \cap \ker \eta \cap T\mathcal{C}_s^{\perp d \eta}) \otimes L^{-1/2} ) \, , $$
we obtain by restriction a complex metric on $H^{\mathbb{C}}$
$$g_{H^{\mathbb{C}}} = d \eta|_{E_m} \otimes \hat{w}_H \, . $$
The quaternionic metric on $M^{\mathbb{C}}$  has the asymptotic development 
\begin{equation*}
g = \frac{1}{\rho^2} g_{-2} + \frac{1}{\rho} g_{-1} + \cdots
\end{equation*}
with 
\begin{equation*}
  g_{-2} = d\rho ^2 + \eta_1^2 + \eta_2^2 + \eta_3^2 \; \mbox{and} \; g_{-1}|_{H^{%
\mathbb{C}}} = g_{H^{\mathbb{C}}} \, .
\end{equation*}
Finally, we put $w_{i,-1} = d \eta_i + \gamma_i$ where $\gamma_i |_{H^{\mathbb{C}}} =
0$.

\subsection{Boundary conditions}

We follow the notations of the previous section and restrict ourselves to
the real slice. We choose an arbitrary complementary $W$ to $H$. Let $(R_1,
R_2, R_3)$ be the dual basis of $(\eta_1, \eta_2, \eta_3)$ on $W$ and let $\tilde{I}_i$ be the almost complex structures on $H$.

The symplectic forms $w_i$ and the metric define almost complex structures $I_i$. Because of
the form of the $w_i$, we have the analytic development 
\begin{equation*}
I_i \partial_{\rho} = I_{i,0} \partial_{\rho} + \rho I_{i,1} \partial_{\rho}
+\cdots = R_i + \psi_i +\cdots 
\end{equation*}
where $\psi_i \in H$ is independent of $\rho$ and if $X \in H$, 
\begin{equation*}
I_i X = I_{i,0} X + \rho I_{i,1} X + \cdots = \tilde{I}_i X + \cdots 
\end{equation*}
We are now in position to show the following

\begin{prop}
The boundary of an AHQK manifold admitting an analytic extension to the boundary is an integrable quaternionic contact structure.
\end{prop}

\begin{proof}
If $X \in H$, one has 
$$w_i( I_j \partial_{\rho} , X ) + w_j ( I_i \partial_{\rho} , X)  = - 2 \delta_{ij} g ( \partial_{\rho} , X ) \, .$$
 The order $-2$ terms do not give anything but from the order $-1$ terms we deduce the equation
$$ w_{i,-1} ( R_j, X ) + w_{j,-1} (R_i, X ) + w_{i,-1} ( \psi_j, X ) + w_{j,-1} ( \psi_i , X )  = - 2 \delta_{ij} g_{-1} ( \partial_{\rho }, X )  $$
so that
$$\begin{array}{rcl} 
d \eta_i (R_j, X ) + d \eta_j( R_i, X)& =&  - \gamma_i ( R_j, X) - \gamma_j (R_i, X ) \\
 & &  - 2 \delta_{ij} g_{-1} ( \partial_{\rho} , X ) + g_{-1} ( \psi_j, \tilde{I}_i X ) + g_{-1} ( \psi_i , \tilde{I}_j X ) \, .
\end{array} $$
The second line gives an element in $S^{3,1} \oplus S^{1,1}$ therefore we need only to look at $\gamma_i$. We will now use the fact that the metric is quaternionic-K\"ahler. Indeed, there exists one forms $\beta_{ij} $ such that the 2-forms $(w_i)$ satisfy
$$d w_i = \sum_j \beta_{ji} \wedge w_j \, , \; \;  \beta_{ji} = - \beta_{ij} \, . $$
 The application $ (\Lambda^1) ^3 \rightarrow \Lambda^3$ 
$$(a_i)_{i=1,2,3} \mapsto \sum_i a_i \wedge w_i$$
 is an injection so that the $\beta_{ij} $ are unique.

We have
$$
\begin{array}{rcl} 
d w_i & =  & - \frac{1}{\rho^3} \sum_{r,s} \varepsilon^{rsi} d\rho \wedge \eta_r \wedge \eta_s \\
& &  + \frac{1}{\rho^2} \left ( d\rho \wedge d\eta_i  + \frac 12 \sum_{r,s} \varepsilon^{rsi} ( d \eta_r \wedge \eta_s - \eta_r \wedge d \eta_s ) \right ) - \frac{1}{\rho^2} d\rho \wedge (d \eta_i + \gamma_i ) +\cdots
\end{array} $$
and then
$$d w_i = - \frac{1}{\rho^3} \sum_{r,s} \varepsilon^{rsi} d \rho \wedge \eta_r \wedge \eta_s + \frac{1}{\rho^2} \left ( \sum_{r,s} \varepsilon^{rsi} d\eta_r \wedge \eta_s - d \rho \wedge \gamma_i \right ) + \cdots \, . $$
We have $\sum_{r,s,p,q} \varepsilon^{irs} \varepsilon^{pqr} \eta_s \wedge \eta_p \wedge \eta_q =0$ so
$$d w_i = \sum_r \left ( \frac {1}{\rho} \sum_s \varepsilon^{irs} \eta_s \right ) \wedge \frac{w_{r,-2} }{\rho^2}  + \frac{1}{\rho^2} \left ( \sum_{r,s} \varepsilon^{rsi} d\eta_r \wedge \eta_s - d \rho \wedge \gamma_i \right ) + \cdots$$
The exterior product of 1-forms with $w_{r,-2}$ is an injection, so $\beta_{ri}$ is of the form
$$\beta_{ri} = \frac{1}{\rho} \beta_{ri,-1} + \beta_{ri, 0} +\cdots \; \mbox{and} \; \; \beta_{ri, -1} = \sum_s \varepsilon^{irs} \eta_s \, . $$

Looking at the order $-2$ terms with respect to $\rho$, one obtains the equations
\begin{equation*} \label{equ} \sum_r \beta_{ri, 0} \wedge w_{r,-2} + \sum_r \beta_{ri, -1} \wedge w_{r,-1}  = \sum_{r,s} \varepsilon^{rsi} d \eta_r \wedge \eta_s - d \rho \wedge \gamma_i \, . 
\end{equation*}
We put $\beta_{ri,0} = \lambda_{ri} d \rho  + \beta_{ri,0}^N $ and $ \gamma_r = d \rho \wedge \gamma_r ^{\rho} + \gamma_r^N$  where $\beta_{ri,0}^N \in T^* N$ and $\gamma_r ^N \in \Lambda^2 T^* N$.

Taking the $d \rho $ component in the previous equation, one gets
\begin{equation*} 
\label{gam} \gamma_ i ^N= \sum_r \eta_r \wedge \beta_{ri,0}^N- \frac 12 \sum_{r,k,s} \lambda_{ri} \varepsilon^{ksi} \eta_k  \wedge \eta_s  + \sum_{r,s} \varepsilon^{irs} \eta_s \wedge \gamma_r^{\rho}  \, . 
\end{equation*}

But then $\gamma_i (R_j ,X ) + \gamma_j (R_i, X) = 0$ and the lemma follows.
\end{proof}

In the two next sections, we will look at integrable quaternionic contact structures in order to show that they are the boundaries of AHQK metrics.

\section{Integrable quaternionic contact structures}

Let $(N, H )$ be a quaternionic contact structure.

In section 2, we computed the derivation of the quaternionic structure in
the direction of $H$. On the other hand, from the identity $d(d\eta_j) (
R_i, X, Y ) = 0$, we obtain 

\begin{eqnarray} 
( \nabla_{R_i} d\eta_j) (X, Y ) & = & \mathfrak{a}( \nabla \alpha_{ij} ) (X,
Y ) - \sum_k \alpha_{ik} \wedge \alpha_{kj} ( X, Y ) \label{der_vert} \\ 
&  & + \sum_k d \eta_j( R_i , R_k ) d \eta_k ( X, Y ) - g ( I_j T_{R_i , X }
+ T_{R_i , I_j X } , Y ) \nonumber \, .
\end{eqnarray}

From now on, we suppose that the quaternionic structure is integrable. We choose a compatible metric $g$ on $H$ and $W =W^g$ the associated complementary vector bundle defining the adapted connection $\nabla$.

\subsection{Torsion}
 The computations of section \ref{horder} give
for any $X \in H$, 
\begin{equation*}
\nabla_X I_j = - \sum_{i=1}^3 \alpha_{ij }(X) I_i \, . 
\end{equation*}

\begin{lm}
\label{S22} Let $(M,H)$ be an integrable quaternionic contact structure. The tensor $T^H$ defined in lemma \ref{1} lies in the component $S^{2,2}$ of $W^* \otimes \mathfrak{so}(H)^{\perp}$.
\end{lm}

\begin{proof}
 By construction, $T^H$ is a section of
$$ \Lambda^2_+ H \otimes \mathfrak{so}(H)^{\perp} = S^{2,0} \oplus S^{4,2} \oplus S^{2,2} \oplus S^{0,2} \, ,$$
so we can put
$$ T_{R_i } = \lambda_i Id + \sum_p I_p A_{pi} \,  $$
with $A_{pi} \in \Gamma ( \Lambda^2_- H)$ ( seen as skew-symmetric endomorphisms ). We apply (\ref{der_vert}) with $i=j$ and obtain $\lambda_i=0$ and $A_{ii} = 0$. Applying one more time (\ref{der_vert}), we see that $A_{pi} $ is equal to the $\Lambda^2_-$ part of $ \mathfrak{a} ( \alpha_{pi}) - \sum_k \alpha_{pk} \wedge \alpha_{ik} $ which is skew-symmetric in $p$ and $i$.

 Writing  $A_i = \frac 12 \sum_{r,s} \varepsilon^{rsi} A_{rs}$, we obtain $T^H$ as the image of $\sum_i I_i \otimes A_i \in \Lambda^2_+ H \otimes \Lambda^2_- H $ by the $SO(4)$-equivariant map $I_i \otimes B \mapsto \frac 12 \sum_j \eta_j \otimes [I_i ,I_j ] B $.
\end{proof}
We are now able to calculate more precisely the vertical derivatives of the
quaternionic structure.

\begin{lm}
\label{der_tot} There exists a function $\lambda$ on $N$ such that 
\begin{equation*}
\nabla_{R_i} I_j = \frac 12 \sum_{k= 1}^3 \left ( d\eta_j ( R_i , R_k ) +
d\eta_i ( R_j , R_k ) - d \eta_k ( R_i , R_j ) \right ) I_k + \lambda [I_i,
I_j ] \, . 
\end{equation*}
\end{lm}

\begin{proof}
  Symmetrizing (\ref{der_vert}) gives
\begin{equation}
\label{sym} \nabla_{R_i} I_j + \nabla_{R_j} I_i = \sum_{k=1} ^3 ( d \eta_j ( R_i , R_k  ) + d \eta_i ( R_j , R_k )  ) I_k \, , 
\end{equation}
In particular, 
$$ \langle \nabla_{R_j} I_j , I_i \rangle = - \langle \nabla_{R_j} I_i , I_j \rangle   = - 2 d\eta_j ( R_i, R_j ) \, , $$
so that we know $\nabla_{R_j} I_i$ except for its component on $[I_i, I_j ]$. We can put 
$$\nabla_{R_i}  I_j = \frac 12 \sum_{k= 1}^3 \left ( d\eta_j ( R_i , R_k  ) + d\eta_i ( R_j , R_k ) - d \eta_k ( R_i , R_j  ) \right ) I_k  + \lambda_{ij}  [I_i, I_j ] \, , $$
with $\lambda_{ii} = 0$. From (\ref{sym}), we have $\lambda_{ij} = \lambda_{ji} $. Moreover, taking for instance $i=1$, $j=2$ and using the skew-symmetry $ \langle \nabla_{R_1} I_2 , I_3\rangle = - \langle \nabla_{R_1} I_3 , I_2 \rangle$, we get $\lambda_{12} = \lambda_{13}$. The other equalities are obtained in the same way. 
\end{proof}

\subsection{The curvature tensor}

We will give some results about the curvature tensor in the $T^{W^g} =0$
case. They will be useful for the twistorial construction.

We are now interested in the curvature $R$ of $\nabla$, and more precisely
in its horizontal part. This is a section $R \in \Gamma ( \Lambda^2 H^*
\otimes \mathfrak{so} (H) )$. The splitting $\Lambda^2 = \Lambda^2_+ \oplus
\Lambda^2_-$ allows us to decompose the curvature in $\Lambda^2_+ \otimes
\Lambda^2_+$, $\Lambda^2_- \otimes \Lambda^2_+$ and $\Lambda^2_- \otimes
\Lambda^2_-$ parts. Looking at its action on $\Lambda^2_+ H$, we have 
\begin{equation*}
\begin{array}{rcl}
R_{X,Y} I_i & = & \nabla_X \nabla_Y I_i - \nabla _Y \nabla_X I_i - \nabla_{
[X,Y]_H} I_i + \sum_j d\eta_j ( X, Y ) \nabla_{R_j} I_i \\ 
& = & \sum_j \left ( - \mathfrak{a} ( \nabla \alpha_{ji} ) + \sum_{k=1} ^3
\alpha_{jk} \wedge \alpha_{ki} \right ) ( X, Y ) I_j + \sum_j d\eta_j ( X, Y
) \nabla_{R_j} I_i
\end{array}
\end{equation*}

\begin{prop}
The $\Lambda^2_+ H \otimes \Lambda^2_+ H$ part of the curvature is scalar.
More precisely, if we denote it by ${\mathcal{S}} \in \Gamma ( End (
\Lambda^2_+ H ) ) $, we obtain with the notations of lemma \ref{der_tot} : 
\begin{equation*}
\mathcal{S} = 2 \lambda Id_{\Lambda^2_+ } \, . 
\end{equation*}
\end{prop}

\begin{proof}
Using Lemma \ref{der_tot} and (\ref{der_vert}), one sees that
$$  ( - \mathfrak{a} ( \nabla \alpha_{ji}  ) + \sum_k \alpha_{jk} \wedge \alpha_{ki}  ) _+ = \frac 12 \sum_k \left ( d \eta_i ( R_j , R_k ) - d \eta_j ( R_i , R_k ) - d \eta_k ( R_i , R_j ) \right ) I_k + \lambda [I_i , I_j ] \, , $$
where the subscript $+$ means the selfdual part. Injecting this in the curvature formula, one easily deduces the proposition.
\end{proof}

We can define a Ricci tensor and a scalar curvature for the partial
curvature $R$. As usual, we put 
\begin{equation*}
\begin{array}{rcl}
Ric(X, Y ) & = & tr_H ( Z\mapsto R_{Z,X} Y ) \\ 
s & = & tr_H ( Ric ) \, ,
\end{array}
\end{equation*}
where the subscript $H$ means that the trace is taken only on $H$. We note $%
Ric_0 $ the trace-free part of the Ricci tensor. In order to obtain the
exact form of the curvature, we use the first Bianchi identity 
\begin{equation*}
\begin{array}{rcl}
R_{X,Y} Z + R_{Y,Z} X + R_{Z,X} Y & = & ( d^{\nabla} T ) _{X,Y,Z}
\end{array}
\end{equation*}

Let $X$, $Y$, $Z$ and $R_i$ be parallel at the point $p$. Since the
horizontal covariant derivatives of $R_i$ and $I_i$ are identical, 
\begin{equation*}
(\nabla T|_H )_p = ( \nabla \sum_i w_i \otimes R_i )_p = 0 \, , 
\end{equation*}
so that at $p$, we have 
\begin{equation*}
\begin{array}{rcl}
R_{X,Y} Z + R_{Y,Z} X + R_{Z,X} Y & = & - T_{[X,Y], Z } - T_{[Y,Z], X } -
T_{[Z,X],Y } \\ 
& = & \sum_{i=1} ^3 ( d \eta_i \wedge T_{R_i} ^H) ( X, Y ,Z )
\end{array}
\end{equation*}

The image by the Bianchi map $b$ of the curvature $R$ lives in the factor $%
S^{2,2} \simeq \Lambda^2_+ H \otimes \Lambda^2_- H$ of $\Lambda^3 H^* \otimes H
\simeq S^{2,0} \oplus S^{0,2} \oplus S^{0,0} \oplus S^{2,2}$.

\begin{prop}
The horizontal part $R \in \Gamma ( \Lambda^2 H^* \otimes \mathfrak{so} (H )
) $ of the curvature tensor seen as an endomorphism of $\Lambda^2 H =
\Lambda^2_+ H \oplus \Lambda^2_- H $ has matrix 
\begin{equation*}
R = \left ( 
\begin{array}{cc}
\frac{s}{12} Id & Ric_0 + B \\ 
{}^t Ric_0 - {}^t B & \frac{s}{12} Id + W^-
\end{array}
\right ) 
\end{equation*}
\end{prop}

\begin{proof}
Recall that the kernel of $b$ is exactly the Riemannian curvature tensors. We have
$$ \begin{array}{rcl} 
S^2 ( \Lambda^2 H^* ) & = & \ker b \oplus \Lambda^4 H^* \\
\Lambda^2 ( \Lambda^2 H^* ) & = & \Lambda^2_+ H^* \oplus \Lambda^2_- H^* \oplus
\Lambda^2_+ H^* \otimes \Lambda^2_- H^ *  
\end{array}$$
We have shown that $b (R) \in S^{2,2}$ so that $R$ is the sum of a Riemannian tensor and an element in the unique irreducible $S^{2,2}$ component which appears in $\Lambda^2 ( \Lambda^2 H^* ) \subset End( \Lambda^2 ( H^* ))$. Moreover, $Ric(B)=0$ if $B \in S^{2,2} \subset \Lambda^2 ( \Lambda^2 H^* )$ so that the Ricci tensor behaves like the Riemannian Ricci tensor, hence is symmetric.
\end{proof}

We show the following lemma which will be useful in the next
section.

\begin{lm}
\label{autoduality} If the vertical torsion vanishes, the curvature $R$ of the adapted connection satisfies the
following equality 
\begin{equation*}
\left [ R_{X, Y} + I R_{X, IY} + I R_{X, IY } - R_{IX, IY} , I \right ] = 0
\, , 
\end{equation*}
for all $X,Y \in H$ and $I \in \Lambda^2_+ H $.
\end{lm}

\begin{proof}
This lemma is well known in the case of anti-selfdual Riemannian curvature in dimension $4$. In our case, it is similar except for the Bianchi part $B$ of the curvature tensor, hence we need only to show that $B$ satisfies the previous equality. We take for instance
$$ B : w \in \mathfrak{so}(H) \mapsto tr ( w K) J - tr ( w J) K \in End ( \mathfrak{so}(H ) ) \, , $$
where $ J \in \Lambda^2_- H$ and $K \in \Lambda^2_+ H$. We must show that
$$ C = \left [ B (w) + I B( Iw) -I B(wI) + B( IwI), I \right ] =0 \, . $$ 

$[J, K ] =0$, hence we get
$$C = \left [ - tr (w J ) K - tr(Iw J ) IK + tr ( w IJ ) IK - tr ( IwIJ ) K, I \right ] \, . $$ 
The result follows then from the two equalities
$$\left \{
\begin{array}{c}
tr ( I w J)  = tr (JI w ) = tr ( IJ w ) = tr ( w IJ ) \\
tr ( Iw IJ ) =tr ( Iw JI ) = - tr ( wJ )
\end{array} \right. $$
\end{proof}

\section{Twistor space}
In this section, we will end the proof of theorem \ref{thm1}.
\subsection{Definitions}

Let $N^7$ be a smooth manifold and $H$ be a quaternionic contact structure
on $N$ with vanishing vertical torsion. Let $g$ be a compatible
Carnot-Carath\'eodory metric, $W$ the adapted complementary distribution and $%
\nabla$ the connection associated to $g$.

Let $\mathcal{T}$ be the set of 2-forms $w \in \Lambda^2_+ H^*$
of norm $\sqrt{2}$. This is a $2$-sphere bundle on $M$ called the twistor space of $(N,H)$. It can be identified with the set of almost complex structures
compatible with $g$ and the orientation. Let $%
\pi$ be the projection $\mathcal{T} \rightarrow N$ and choose a
local quaternionic structure $(I_1, I_2, I_3 )$ associated to the $1$-forms $%
(\eta_1, \eta_2 ,\eta_3 )$. At a point $I= x_1 I_1 + x_2 I_2 + x_3 I_3 $, we
put 
\begin{equation*}
\eta^r = x_1 \pi^* \eta_1 + x_2 \pi^* \eta_2 + x_3 \pi^* \eta_3 \, . 
\end{equation*}
It is a well defined 1-form on $\mathcal{T}$ not depending on the choice of $%
SO_3$-trivialization $(I_1, I_2, I_3 )$.

Using the connection $\nabla$, we split the tangent bundle of $\mathcal{T}
$ at $I \in \pi^{-1} (s) $ for $s \in N$:
\begin{equation*}
T_I \mathcal{T} = T_I \mathcal{T}_s \oplus \pi^* T_s N \, . 
\end{equation*}
Here $\mathcal{T}_s$ is the fiber above $s$ of
the fibration $\pi$. We call $Hor_I \mathcal{T} \simeq T_s N = W_s \oplus
H_s$ the horizontal space. Let $(R_1, R_2, R_3)$ be the dual basis of $(\eta_1,
\eta_2 , \eta_3 )$ on $W$. At $I_1$, we have an almost complex structure $J
$ on $\ker \eta_r \simeq \ker \eta_1 \oplus T_I \mathcal{T}_s$ satisfying

\begin{itemize}
\item on $\ker \eta_1$, the almost complex structure satisfies $J = I_1$ after extending $I_1$ to all $\ker \eta_1$
by $I_1 R_2 = R_3$ and $I_1 R_3 = - R_2$;
\item on $T_{I_1} \mathcal{T}_s$ $J $ is the natural complex structure given by
the metric and the orientation on the sphere $\mathcal{T}_s$.
\end{itemize}

\begin{prop}
Let $H$ be an integrable quaternionic contact structure on a 7-dimensional manifold $N$. The almost complex structure $J$ defined on
the kernel of $\eta^r$ is independent of the choice of compatible metric $g$ on $H$.
\end{prop}

\begin{proof}
Let $\eta'= f^2 \eta $ be a conformal change. The distribution $\ker \eta_r$ on the twistor space is left unchanged. The conformal change gives a new complementary $W^{f^2 g }$ spanned by $(R_1', R_2', R_3 ')$ and a new connection $\nabla' = \nabla + a $. The distribution $Hor'_I \mathcal{T}$ is the horizontal subspace on $\mathcal{T}$ corresponding to $\nabla'$, and $J'$ is the corresponding almost complex structure.

 The vertical part of $J$ is left unchanged. 

At $I_1 \in \mathcal{T}_s$, we take $U \in \ker \eta^r$, horizontal for the connection $\nabla$, and $X$ its projection on $N$. In the decomposition $T_{I_1} \mathcal{T} =Hor_{I_1} \mathcal{T} \oplus T_{I_1}  \mathcal{T}_s$, we have $U = ( X, 0)$ and $JU = ( I_1X, 0 )$. On the other hand, in the decomposition $ T_{I_1}   \mathcal{T} = Hor_{I_1} ' \mathcal{T}  \oplus T_{I_1}  \mathcal{T}_s$, we have $ U = ( X, - a_X I_1) $, $JU = ( I_1 X, - a_{I_1 X} I_1  )$ and $J'U = ( I_1 X, - \frac 12 [ I_1  , a_X I_1 ] ) $ thus $J$ and $J'$ coincide iff $ a_{I_1 X} I_1  = \frac 12 [I_1 , a_X I_1 ]$ for all $X \in \ker \eta_1$. 

One has $ a \in \Omega^1 ( \mathbb{R} \oplus \mathfrak{so}(H))$, and we decompose the $\mathfrak{so}(H)$-part in selfdual and anti-selfdual part that we write respectively $a^+$ and $a^-$. From \ref{changeconnection}, one gets 
$$a_X^+ = \sum_j \langle I_j \theta^{\sharp} , X \rangle I_j \, . $$
and verifies easily that $a_{I_1 X} I_1 = \frac 12 [I_1, a_X I_1 ]$.

The quaternionic contact structure has no vertical torsion, so that we get from lemma \ref{changeconnection}
$$ a_{R_i} = \theta(R_i ) + 2 |\theta^{\sharp} |^2 I_i + 2 \theta^{\sharp} \wedge I_i \theta^{\sharp} + 2 ( I_i \nabla \theta^{\sharp} )^{\mathfrak{so}(H)} \, . $$
Taking the selfdual part gives
$$a_{R_i} ^+ = 3 | \theta^{\sharp} |^2 I_i- \frac 12 \sum_k tr ( I_k I_i \nabla \theta^{\sharp} ) I_k \, ,  $$
and computing $a_{R_3} I_1$ and $a_{R_2} I_1$ one obtains $a_{R_3} I_1 = \frac 12 [I_1, a_{R_2} I_1 ] $. 
\end{proof}

\subsection{Integrability of the twistor space}

This section is devoted to the proof of the following theorem :

\begin{thm} \label{CR}
Let $H$ be a quaternionic contact structure with vanishing vertical torsion and $J$
be the almost complex structure on the kernel of $\eta^r$ on the twistor space.
Then

\begin{itemize}
\item  $J$ is adapted to the symplectic form $d \eta^r$ on $\ker \eta^r$ and gives a metric of signature $(6,2)$.

\item  $J$ is integrable.
\end{itemize}
\end{thm}

\begin{proof}
 The first point is similar to \cite{Biq00} and 
$$ d \eta^r ( \cdot, J \cdot ) = g_H + d \eta_1 ( R_2 ,R_3 ) ( \eta_2 ^2 + \eta_3 ^2 ) + \eta_3 \odot d x_2 - \eta_2  \odot d x_3 \, , $$
where $\alpha \odot \beta = \frac 12 ( \alpha \otimes \beta + \beta \otimes \alpha )$ is the symmetric product. This is the metric of signature $(6,2)$.

\bigskip

We must now verify the integrability of $J$. This is given by the vanishing of the Nijenhuis tensor
$$N(X, Y)  = [X,Y ] + J [X, JY ] + J [JX, Y ]- [ JX, JY ] \, . $$ 

 If $X$ and $Y$ are vertical, it follows from the fact that $J$ is the complex structure of the 2-sphere which is integrable, and if $X$ is horizontal and $Y$ vertical this is similar to the proof of 14.68 in \cite{Bes87}.

Assume now that $X$ and $Y$ are horizontal. In this case the vertical part and the horizontal part of $N(X,Y)$ at $I \in \mathcal{T}$ are given by
$$\begin{array}{rcl}
(N(X,Y))_{Hor} & =& T_{X,Y} + I T_{X,IY} +I T_{IX,Y} - T_{IX, IY} \\
(N(X,Y))_{Ver} & =& \left [ R_{X, Y} + I R_{X, IY} + I R_{IX, Y} - R_{IX, IY } , I \right ] 
\end{array}$$

\bigskip

We look first at the horizontal part. If $X, Y \in H$, then $T_{X,Y} = \sum_i d \eta_i (X, Y) R_i$ and we deduce easily that $(N(X,Y))_{Hor} = 0$. If $X= R_2$ and $Y= R_3$ at $I=I_1$, then $(N(X,Y))_{Hor} = T_{R_2, R_3} - T_{R_3, -R_2} =0$ so that the only no-trivial case at $I_1$ is $X \in H$ and $Y =R_2$. Following the notations of \ref{S22}, the $W$-part of the torsion $T_{R,X}$ vanishes and the $H$-part is $T_{R_i,X} = \sum_{p} I_p A_{pi} X$ where $A_{pi} = - A_{ip} \in \Lambda^2_- H$. Therefore, we have
$$\begin{array}{rcl}
(N(X,R_2))_{Hor} & = & - \sum_p I_p A_{p2} X -  I_1 \sum_p I_p A_{p3} X - I_1 \sum_p I_p A_{p2} I_1 X + \sum_p I_p A_{p3} I_1 X \\
 & = & -I_3 A_{32} X -I_3 A_{23} X -I_1 I_3 A_{32} I_1 X +I_2 A_{23} I_1 \,  .
\end{array}$$
The $A_{ij}$ and $I_k$ commute and the skew-symmetry $A_{23} = - A_{32}$ gives the vanishing of $( N_{X,R_2})_{Hor}$. 

\bigskip

We show now the vanishing of the vertical part. If $X, Y \in H$, this is just lemma \ref{autoduality}. It remains to show that for $X \in H$, 
$$c_{R_2, X, Y} = \left [ R_{R_2, X}  + I_1 R_{R_3, X} + I_1 R_{R_2, I_1 X } - R_{R_3, I_1 X} , I_1 \right ]Y = 0 \, . $$
We put $I_1R_1 =0$ and $I_1 R_2 = R_3$, in order to have $c_{X,Y, Z}$ defined for all $X$, $Y$ and $Z$. Because we have the same identities on the torsion, the computation is very similar to \cite[Lemma II.5.3]{Biq00} and one gets
$$ c_{R_2, X ,Y } + c_{Y, R_2, X} = 0 \, , \; \; \forall X ,Y  \in H \, . $$ 

$c_{R_2, X, Y} $ is in the subspace spanned by $I_2 Y $ and $I_3 Y$ therefore if the $\mathbb{C}$-subspaces spanned by $Y$ and $X$ for the almost complex structure $I_1$ are transverses, then $c_{R_2, X, Y} =0$. We deduce that $c_{R_2, X, Y} =$ in all cases.

\end{proof}

\subsection{Proof of theorem 1.1}

We have shown that any integrable quaternionic contact structure $H$ admits a twistor space $\mathcal{T}$ which is CR-integrable. This is sufficient to apply the results of Biquard \cite{Biq00} which give the theorem \ref{thm1} (see part III for the twistorial construction). The idea of the proof is to construct a bigger twistor space $\mathcal{N}$ which has the properties of section \ref{3.1} and such that $\mathcal{T}$ is a real hypersurface of $\mathcal{N}$.
 With the notations of \ref{3.1}, the AHQK metric is $g = w_E \otimes w_H$ and is quaternionic-K\"ahler, \cite{LeB89}.

 The corollary \ref{cor} follows immediately from our theorem \ref{CR} and the theorem 0.4 of \cite{Biq02}.

\section{Deformations of the $7$-sphere}

Hereafter, we assume that $N= \mathbb{S}^7$ is the $7$-sphere in $\mathbb{H}^2$ where $\mathbb{H}^2$ is an $\mathbb{H}$-vector space with $\mathbb{H}$ acting on right. Let $ \langle \cdot , \cdot \rangle $ be the canonical metric on $\mathbb{H}^2 \simeq \mathbb{R}^8$. Recall that we have a quaternionic contact structure on $\mathbb{S}^7$ given by $H_x^c = ( x \mathbb{H} ) ^{\perp}$ for $x \in \mathbb{S}^7$. The restriction to $H^{can}$ of the round metric on $\mathbb{S}^7$ defines an adapted metric $g_0$. The adapted complementary is $W_x = x Im \mathbb{H}$ and is spanned by $R_1(x)  = xi$, $R_2 (x) = xj$ and $R_3 (x) = xk$.

$H^{can}$ is a connection on the principal $Sp(1)$-bundle $ \mathbb{S}^7 \rightarrow \mathbb{S}^4$ (Hopf-bundle). We call $\eta \in \Omega^1 ( \mathbb{S}^7) \otimes \mathfrak{sp}(1) \simeq \Omega^1 (W ) $ its connection form. Let us write it $\eta = \sum_i \eta_i$ or $\sum_i \eta_i \otimes R_i$. One has $d \eta_i ( W, H^{can} ) = 0$ so that the the torsions $T^W = - \frac 12 \sum_{i,j}( \alpha_{ij} + \alpha_{ji} ) \otimes w_i^* \otimes w_j$ and $T^H$ vanish.

 Let $\nu$ be the canonical volume form of $\mathbb{S}^7$ that we decompose as $\nu = \nu^c \wedge \eta_1 \wedge \eta_2 \wedge \eta_3$ so that $\nu^c|_H$ is a volume form on $H^{can}$.

In this section, we compute the complex of integrable infinitesimal deformations of $H^{can}$.

\subsection{Deformation of the integrability condition} \label{int}
A deformation of $H^{can}$ is given by a $1$-form $\theta$ with values in $W$ which vanishes on $W$, or equivalently by a section of $End(H^{can}, W)$. The link between the new distribution and $\theta$ is given by 
$$ H_{\theta} = \{X - \theta(X) , X \in H^{can} \}= \ker ( \eta + \theta ) \, . $$ 

 Assume now that $\theta^t$ is a $1$-parameter family of such $1$-forms, each giving a vertical torsion free distribution denoted by $H_t = \ker ( \eta + \theta ^t ) $. For small $t$,  the forms $ d ( \eta_i + \theta^t_i ) |_{H_t} \in \Gamma ( \Lambda^2  H_t ^* )$ span a space of selfdual $2$-forms on $H_t$ with respect to a metric $g_t$ on $H_t$. We choose $g_t$ such that $g_0$ is the restriction of the round metric on $H^{can}$.

In order to write the condition on the torsion, one has to take an orthonormal basis of $\Lambda^2_+ H_t^*$. We identify the functions and the $4$-forms on $H^{can}$ using $\nu^c$. We search $a^t : \mathbb{S}^7 \rightarrow GL(3, \mathbb{R})$ such that $a^0 =Id$ and 
$$ \left [ a^t \cdot d ( \eta + \theta^t ) \right ]_i \wedge 
   \left [ a^t \cdot d ( \eta + \theta^t ) \right ]_j \wedge
   \left [ a^t \cdot ( \eta + \theta^t ) \right ]_1 \wedge
   \left [ a^t \cdot ( \eta + \theta^t ) \right ]_2 \wedge
   \left [ a^t \cdot ( \eta + \theta^t ) \right ]_3 = 2 \delta_{ij} \nu \, . $$

Setting $ \dot{\psi} = \frac{d \psi^t}{d t} |_{t=0}$, one obtains
\begin{equation}
 \dot{a}_{ij} + \dot{a}_{ji} + ( d \dot{ \theta}_i \wedge d \eta_j + d \dot{\theta}_j \wedge d \eta_i )|_{H^{can}} + tr ( \dot{a} ) = 0 
\end{equation}

\begin{rem}
We used the fact that $\alpha_{ij} = 0$. In general, one has
\begin{equation*}
 \dot{a}_{ij} + \dot{a}_{ji} + ( d \dot{ \theta}_i \wedge d \eta_j + d \dot{\theta}_j \wedge d \eta_i )|_{H} + \sum_k ( \alpha_{ki} \wedge \dot{\theta}_k \wedge d \eta_j  + \alpha_{kj}\wedge \dot{\theta}_k \wedge d \eta_i ) |_{H^{can}} + tr ( \dot{a} ) = 0 \, .
\end{equation*}

\end{rem}

 We put $\beta^t = a^t \cdot (\eta + \theta^t )$ with dual basis $(R_1^t, R_2^t, R_3^t )$ on $W$. Our choice of $a^t$ ensures that we obtain an orthonormal direct basis in $\Lambda^2_+ H_t$ for the metric $g_t$. Let $I_i^t$ be the associated quaternionic structure. By \ref{ifint}, the deformation preserves the integrability iff there exist $\gamma_i^t$ such that for $X \in H^{can}$, 
$$i_{R_i^t}\beta_j ^t (X - \theta^t (X) ) + i_{R_j^t} \beta_i^t ( X - \theta^t (X) ) = \gamma_i ^t \circ I_j^t ( X - \theta^t (X) ) + \gamma_j ^t \circ I_i^t ( X- \theta^t (X)) \, . $$   
The $\gamma_i^0$ vanish so that one obtains the following lemma.

\begin{lm}
If $\theta^t$ is a $1$-parameter smooth deformation of the quaternionic contact structure on $\mathbb{S}^7$ which preserves the integrability, we have
 $$ \mathcal{A}_0(\dot{\theta} ) = - d ( \dot{a}_{ij} + \dot{a}_{ji}  )|_{H^{can}}  + (i_{R_i} d \dot{\theta}_j + i_{R_j} d \dot{\theta}_i ) |_{H^{can}} \in S^{3,1}\oplus S^{1,1} \, , $$
where
$$ \dot{a}_{ij} + \dot{a}_{ji} + ( d \dot{ \theta}_i \wedge d \eta_j + d \dot{\theta}_j \wedge d \eta_i )|_{H^{can}} + tr (\dot{a})=0 \, .$$
\end{lm}

\begin{rem}
The statement has exactly the same form if one deforms Einstein selfdual Levi-Civita connections with non-zero scalar curvature ( which give $3$-Sasakian manifolds and so integrable quaternionic contact structures, see \cite{Kon} ).
\end{rem}

 The composition of $\mathcal{A}_0$ with the projection on $S^{5,1} $ gives a differential operator $\mathcal{A}: \Gamma ( (H^{can})^* \otimes W ) \rightarrow \Gamma ( S^{5,1} ) $. Its kernel gives the infinitesimal deformations of $H^{can}$ preserving the integrability. This kernel contains the image of the infinitesimal diffeomorphisms through 
$$ \begin{array}{rcl}
 \mathcal{D} : \Gamma ( T \mathbb{S}^7 ) & \rightarrow & \Gamma ( (H^{can})^* \otimes W ) \\
 \zeta & \mapsto & \{ X \in H \mapsto X. \eta( \zeta ) + d \eta ( \zeta , X ) \} 
\end{array} $$


\subsection{A Bianchi identity} Because of the dimensions of the different vector bundles, the previous complex cannot be elliptic, even in the direction of $H^{can}$. We will show now a Bianchi identity.

\begin{lm} \label{bianchi}
 Let $(M^7,H,g)$ be a quaternionic contact structure where $g$ is a particular choice of Carnot-Carath\'eodory metric. Let $W$ be the adapted complementary and $\nabla$ be the corresponding adapted connection. The vertical torsion $T^W$ of $H$ is a section of $S^{5,1} \subset H^* \otimes S^{4,0}$. Let $\mathcal{B}_H$ be the composition of $d^{\nabla} : \Gamma( H^* \otimes S^{4,0} ) \rightarrow \Gamma ( \Lambda^2 H^* \otimes S^{4,0} ) $ with the projection on $S^{6,0}$. Then we have
$$ \mathcal{B}_H ( T^W ) =0 \, . $$ 
\end{lm}

\begin{rem}
Here is a small abuse of notation. Indeed $d^{\nabla}$ can be applied only on true $1$-forms with values in a vector bundle. Nevertheless we can give the following meaning to $d^{\nabla}$ : a section $\sigma$ of $H^* \otimes E$ is extended in a true $1$-form vanishing on $W$ and we use then the vanishing $(T_{X,Y})_H =0$ in order to obtain
$$\begin{array}{rcl}
( d^{\nabla} \sigma ) (X,Y) & = &\nabla_X \sigma_Y - \nabla_Y \sigma_X - \sigma_{T_{X,Y} } \\
 & = & ( \nabla_X \sigma) _Y - ( \nabla_Y \sigma ) _X \, , 
\end{array} $$
for vector fields $X,Y \in H$.
This kind of equalities will be used throughout the proof for every elements of $ \Gamma ( H^* \otimes E)$ and every vector bundle $E$.
 
\end{rem}

\begin{proof}
 Let $(I_1, I_2, I_3) $ be a local direct orthonormal basis of $\Lambda^2_+ H^*$ corresponding to local $1$-forms $(\eta_1, \eta_2 , \eta_3 ) $ defining the contact structure. Denote by $(R_1, R_2, R_3)$ the corresponding dual basis on $W$. The first Bianchi identity gives
$$\mathfrak{S}_{X,Y,R_i} \left ( R_{X,Y} R_i -T_{T_{X,Y} , R_i } - ( \nabla_X T ) _{Y,R_i }  \right ) = 0 \, , $$
where $X$ and $Y$ are two vector fields in $H$. Taking the $W$-part, we obtain
\begin{equation*}
R_{X,Y} R_i = (T_{T_{R_i,X} , Y })_W + ( T_{T_{Y, R_i}, X } )_W +( ( \nabla_X T )_{Y,R_i})_{W} + (( \nabla_{R_i} T)_{X,Y} )_W + (( \nabla_Y T)_{R_i, X } )_W \, .      
\end{equation*}

 We calculate first $A_1 ( X, Y ,R_i ) =  (T_{T_{R_i,X} , Y })_W + ( T_{T_{Y, R_i}, X } )_W $. One has
$$A_1 ( X, Y , R_i ) = - T^W_{T^W_{X}(R_i)} (Y) + \sum_j \langle I_j  T^H_{R_i} (X), Y \rangle R_j  + T^W_{T^W_ {X} (R_i) } (X ) -\sum_j \langle I_j T^H_{R_i} (Y) , X \rangle R_j \, . $$
We put $a_{ij} =\frac 12 ( \alpha_{ij} + \alpha_{ji} )$. With this notations, we get
\begin{equation}
\label{b1}
 A_1 ( X, Y , R_i ) = - \sum_{k,j=1} ^3 a_{ji} \wedge a_{kj} (X, Y ) R_k  + \sum_{k=1} ^3 \langle ( I_k T^H_{R_i} + T^H_{R_i} I_k ) (X) , Y \rangle R_k \, . 
\end{equation}

We suppose now that $p \in M$, and that $X$, $Y$ and $R_i$ are parallel at $p$. In particular, we have at $\alpha_{ij} = \alpha_{ji}$ at $p$. Then at $p$, 
$$( ( \nabla_X T )_{Y,R_i} + ( \nabla_{R_i} T)_{X,Y}  + ( \nabla_Y T)_{R_i, X } )_W = - \sum_{k} ( d^\nabla a_{ki} ) (X, Y) R_k - \sum_k ( \nabla_{R_i} d \eta_k ) (X, Y ) R_k \, ,$$
so that we obtain
$$\begin{array}{rcl}
R_{X,Y} R_i & = & - \sum_{k,j=1} ^3 a_{ji} \wedge a_{kj} (X, Y ) R_k  + \sum_{k=1} ^3 \langle ( I_k T^H_{R_i} + T^H_{R_i} I_k ) (X) , Y \rangle R_k \\
 & &  - \sum_{k} ( d^\nabla a_{ki} ) (X, Y) R_k - \sum_k ( \nabla_{R_i} d \eta_k ) (X, Y ) R_k \, .
\end{array} $$

From the obvious equation $ \langle R_{X, Y} R_i , R_k \rangle + \langle R_{X,Y} R_k , R_i \rangle =0 $, we deduce that
$$\begin{array}{rcl} 
 2 d ^{\nabla} a_{ki} (X, Y ) & =  &  \langle ( I_k T^H_{R_i} + T^H_{R_i} I_k ) (X) , Y \rangle +  \langle ( I_i T^H_{R_k} + T^H_{R_k} I_i ) (X) , Y \rangle
\\
 & & - \nabla_{R_i} d \eta_k (X,Y ) - \nabla_{R_k} d \eta_i (X, Y ) \, .
\end{array}  $$

Remark that (\ref{der_vert}) is true even if $T^W$ does not vanish. At $p$, it gives

$$ \begin{array}{rcl}

 4 d^{\nabla} a_{ki}  & = &  \sum_j ( d \eta_k ( R_j , R_i ) + d \eta_i (R_j , R_k )) \langle I_j \cdot , \cdot \rangle \\
& &  + 2\langle ( I_k T^H_{R_i} + T^H_{R_i} I_k ) \cdot , \cdot \rangle +  \langle ( I_i T^H_{R_k} + T^H_{R_k} I_i )  \cdot , \cdot \rangle 
\end{array} $$

This is a $2$-form and taking the selfdual part, we obtain
$$4 ( d^{\nabla} a_{ki} )_+ =  \sum_j ( d \eta_k ( R_j , R_i ) + d \eta_i (R_j , R_k )) \langle I_j \cdot , \cdot \rangle- 2\, tr ( T_{R_i}^H) \langle I_k \cdot , \cdot \rangle    - 2\, tr ( T^H_{R_k} )  \langle I_i \cdot , \cdot \rangle  \, . $$

This is an element of $S^{2,0} \otimes (S^{4,0} \oplus S^{0,0}) \subset (S^{2,0} ) ^3 $. We take the projection in $Sym^3 ( S^{2,0} ) \simeq S^{6,0} \oplus S^{2,0} $ and then the $S^{6,0}$-part to obtain the lemma. 

\end{proof}

\subsection{The complex of infinitesimal deformations } We take the infinitesimal part of the previous equation and obtain the complex of infinitesimal deformations of the $7$-sphere
$$\begin{array}{rcc}
(\mathcal{C}_0 ) &  \Gamma (T \mathbb{S}^7 ) \stackrel{\mathcal{D} }{\rightarrow} \Gamma(H^* \otimes W ) \stackrel{\mathcal{A}}{\rightarrow} \Gamma ( S^{5,1} ) \stackrel{\mathcal{B}_c}{\rightarrow } \Gamma ( S^{6,0} ) \, . &
\end{array}$$
 Here $\mathcal{B}_c$ means the Bianchi operator on $H^{can}$.

We have the decomposition $\Gamma(T \mathbb{S}^7 ) =\Gamma(W) \oplus \Gamma (H^{can})$ and on the other hand $\Gamma( (H^{can})^* \otimes W ) = \Gamma ( S^{3,1}) \oplus \Gamma (S^{1,1})$ with the property that $\mathcal{A}(\Gamma(S^{1,1} )) =0$. The restriction of $\mathcal{D}$ to $\Gamma ( H^{can}) $ is an isomorphism $ \Gamma (H^{can}) \rightarrow \Gamma(S^{1,1} )$ so that if $\tilde{\mathcal{D}}$ is the composition of $\mathcal{D}$ restricted to $\Gamma ( W)$ with the projection on $\Gamma (S^{3,1})$, we obtain an isomorphism
$$ \frac{\ker \mathcal{A}}{ \mathcal{D}( \Gamma(T \mathbb{S}^7 ))} \simeq \frac{\ker \mathcal{A} \cap \Gamma ( S^{3,1} )}  {\tilde{\mathcal{D}} ( \Gamma ( W))} \, . $$ 
In other words, we can compute the first homology group of the complex
$$\begin{array}{lcc}
( \mathcal{C} )  & \Gamma ( W) \stackrel{\tilde{\mathcal{D}} }{\rightarrow} \Gamma(S^{3,1} ) \stackrel{\mathcal{A}}{\rightarrow} \Gamma ( S^{5,1} ) \stackrel{\mathcal{B}_c}{\rightarrow} \Gamma ( S^{6,0} )  \, .  &
\end{array}$$

\begin{rem}
This complex is not elliptic. Nevertheless a straightforward computation shows that $(\mathcal{C})$ is elliptic in the direction of $H^{can}$. This was not the case of $(\mathcal{C}_0)$.
\end{rem}
\begin{lm} If $\xi= \xi_{H^{can}} + \xi_W  \in T_x \mathbb{S}^7$, the principal symbols $\sigma_{\xi}$ of the previous differential operators satisfy :
\begin{itemize}
\item If $\xi_{H^{can}}= 0$, then $\ker \sigma_{\xi} ( \tilde{\mathcal{D}} ) = W_x $, or else if $\xi_{H^{can}} \neq 0$, then $\ker \sigma_{\xi} ( \tilde{\mathcal{D}} ) = \{0 \}$.
\item If $\xi_{H^{can}} =0$, then $\ker \sigma_{\xi} ( \mathcal{A} ) = S^{3,1}_x$, or else if $\xi_{H^{can}} \neq 0$, then $\ker \sigma_{\xi} ( \mathcal{A} )  = \mbox{Im} \, \sigma_{\xi} ( \tilde{\mathcal{D}} )$.
\item If $\xi_{H^{can}} =0$, then $\ker \sigma_{\xi} ( \mathcal{B} ) = S^{5,1}_x$, or else if $\xi_{H^{can}} \neq 0$, then $\ker \sigma_{\xi} ( \mathcal{B}_c ) = \mbox{Im} \, \sigma_{\xi} ( \mathcal{A}) $.
\end{itemize}
\end{lm}

\section{ $Sp(1)$-invariant deformations of the $7$-sphere}

We have seen in the previous section that infinitesimal deformations of the standard quaternionic contact structure on $\mathbb{S}^7$ are parametrized by the first cohomology group of the complex $(\mathcal{C})$. This complex is not elliptic and even not hypoelliptic. Indeed, \cite{LeB91} ensures the existence of an infinite dimensional moduli space of integrable quaternionic contact structures on $\mathbb{S}^7$.

 In order to obtain an elliptic complex, we will look at quaternionic contact structures on $\mathbb{S}^7$ admitting a free $Sp(1)$-action. Here, $Sp(1)$ is viewed as the group of unitary quaternions. There is a canonical action of $Sp(1)$ on $\mathbb{S}^7$ given by the  diagonal action of $Sp(1)$ on $\mathbb{S}^7 \subset \mathbb{H}^2$. The quotient is the $4$-sphere and the projection $\mathbb{S}^7 \rightarrow \mathbb{S}^4$ is the Hopf projection. Smooth deformations of this $Sp(1)$-action on $\mathbb{S}^7$ are always diffeomorphic to the canonical one. Therefore, we fix the $Sp(1)$-action to be the canonical one. 

\subsection{$G$-invariant structures } \label{G-invariant} In this section, we do some general remarks about quaternionic contact structures $H$ invariant under a free smooth $G$-action, where $G= SO(3)$ or $G= Sp(1)$. Let $(N,H)$ be such a quaternionic contact structure. The action must be transverse to the contact distribution so that $H$ is a connection on a $G$-principal bundle $ N  \stackrel{\pi} {\rightarrow} B$. Let $(\eta_1, \eta_2 , \eta_3 )$ be the connection form of $H$ with values in $\mathfrak{sp}(1)$. The symplectic forms $d \eta_i |_H$ define a unique adapted conformal class of metrics $[g]$ on $H$. Because of the $G$-invariance, the conformal class $[g]$ can be pushed down on $B$ and gives a conformal class of Riemannian metrics $[g]$ on $B$. Let $E = M \times_{Ad} \mathfrak{g}$ be the adjoint bundle. The connection $H$ gives a covariant derivative $\nabla^E$ on $E$, with curvature $R^E$. By definition of $[g]$, the curvature $R^E$ gives an isomorphism
$$ \begin{array}{rcl}
R^E : \Lambda^2_+ T B & \rightarrow & E \\
  \zeta & \mapsto & R_{\zeta} \in \mathfrak{so}(E) \simeq E 
\end{array} $$


Let $D$ be a linear connection, preserving the conformal class. Every choice of $D$ is available, but in general one chooses a metric $g$ in the conformal class and the corresponding Levi-Civita connection. 

The tensor $(R^E)^{-1} \nabla^{D,E} R^E$ is a section of $T^* B \otimes End(\Lambda^2_+ B) $ and taking the symmetric part of $End(\Lambda^2_+ B )$ with respect to any choice of metric in the conformal class $[g]$, it gives a tensor $T$ in $ \Gamma ( S^{1,1} \otimes ( S^{4,0} \oplus S^{0,0} ))$. The $S^{5,1}$ part of $T$ is the vertical torsion of the quaternionic contact structure $H$. We put $T\! or ( H ) =T $. 

\subsection{Infinitesimal $Sp(1)$-invariant deformations of $\mathbb{S}^7$ } We now come back to the deformations of the canonical quaternionic contact structure on $\mathbb{S}^7$. Let $\mathcal{H}$ be the set of $Sp(1)$-invariant quaternionic contact structures on the Hopf bundle $\mathbb{S}^7 \rightarrow \mathbb{S}^4$ and $\mathcal{G}$ be the group of diffeomorphisms of $\mathbb{S}^7$ commuting with the $Sp(1)$ action. Let $\nabla$ be the Levi-Civita connection of the round metric of $\mathbb{S}^4$. In the $Sp(1)$-invariant case, the complex $(\mathcal{C})$ can be written on the basis $\mathbb{S}^4$ in the following way :

\begin{lm} \label{sp1complex}
The complex $(\mathcal{C})$ applied to $Sp(1)$-invariant deformations on the Hopf bundle $\mathbb{S}^7 \rightarrow \mathbb{S}^4$ can be written on the basis as
$$\begin{array}{lcl}
 (\mathcal{C}) & \Gamma(S^{2,0}) \stackrel{\tilde{\mathcal{D}}}{\rightarrow} \Gamma (S^{3,1} ) \stackrel{\mathcal{A}}{\rightarrow } \Gamma ( S^{5,1} ) \stackrel{\mathcal{B}_c}{\rightarrow} \Gamma( S^{6,0}) & 
\end{array} $$ 
where $\tilde{\mathcal{D}} = p^{3,1} \nabla$, $\mathcal{A} = p^{5,1} \nabla^2$ and $\mathcal{B}_c = p^{6,0} \nabla $. The homology groups $H^0$, $H^1$, $H^2$ and $H^3$ of $(\mathcal{C})$ have dimensions $10$, $35$, $0$ and $0$ respectively.
\end{lm}

\begin{proof}
The operator $\mathcal{A}$ is the composition of $\mathcal{A}_1 = p^{4,0} \nabla$ and $\mathcal{A}_2 = p^{5,1} \nabla$ so that the previous complex splits into
$$ \begin{array}{rcl}
\label{c1} (\mathcal{C}_1) &   \Gamma(S^{2,0}) \stackrel{\tilde{\mathcal{D}}}{\rightarrow} \Gamma( S^{3,1} ) 
\stackrel{\mathcal{A}_1}{\rightarrow} \Gamma(S^{4,0}) \, , &  \\

\label{c2}  (\mathcal{C}_2) & \Gamma(S^{4,0}) \stackrel{\mathcal{A}_2} {\rightarrow} \Gamma (S^{5,1} ) \stackrel{\mathcal{B}_c}{\rightarrow} \Gamma( S^{6,0}) \, . & 
\end{array}$$

One recognize in $(\mathcal{C}_1)$ the complex of deformations of anti-selfdual metrics. This complex is well-known and one can show that $\mathcal{A}_1$ gives an isomorphism between $\ker \tilde{\mathcal{D}}^*$ and $\Gamma(S^{4,0} )$ (see for instance the proof of \cite[theorem 13.30, p. 376]{Bes87} ). Therefore the kernel of $\mathcal{A} \oplus \mathcal{D}^* $ can be identified with the kernel of $\mathcal{A}_2$, and we are reduced to the study of $(\mathcal{C}_2)$.

First we give some Weitzenb\"ock formula. Let $\mathcal{D}^{\nabla}$ be the Dirac operator on $S \otimes E$ where $S =S^{1,0} \oplus S^{0,1}$ is the spinor bundle and $E$ can be any $S^{n,m}$. The Dirac operator is the composition of the connection and the Clifford multiplication. The Clifford multiplication is a morphism of representation on $Spin(4)$-modules so is the identity on each irreducible component of $S \otimes E$, up to a multiplicative constant. If $E=S^{5,0}$, we see for instance that $\mathcal{D}^{\nabla} = b \mathcal{B}_c \oplus a \mathcal{A}_2^* $ for some constants $a$ and $b$. The Weitzenb\"ock formula is 
$$\mathcal{D}^{\nabla}( \phi \otimes s ) = \nabla^* \nabla ( \phi \otimes s ) + \frac{s}{4} \phi \otimes s + \sum_{e_i, e_j } e_i \cdot e_j \phi \otimes R_{e_i,e_j}^{\nabla} s \, , $$ 
and in our case, the curvature $R^{\nabla}$ is scalar so that the last term in the previous equality is a combination of Casimir operators ( see \cite[p. 376]{Bes87}. One obtains finally 
$$(\mathcal{D}^{\nabla} ) ^2 = \nabla^* \nabla + \frac{s}{4} \, ,$$
and so $\ker( \mathcal{B} \oplus \mathcal{A}_2^* ) = 0$, that is to say the complex $(\mathcal{C})$ has no second homology group. 

In the same way, regarding $\mathcal{B}_c^* : \Gamma(S^{6,0} ) \rightarrow \Gamma ( S^{5,1} ) $, it appears to be the Dirac operator on $S^{6,0} \subset S^{1,0} \otimes S^{5,0} $ ( up to a multiplicative constant ). One can show that 
$$ (( \mathcal{D}^{\nabla} )^2- \nabla^* \nabla ) |_{S^{6,0} } = \frac{2s}{3} \, , $$
which gives $\ker ( \mathcal{B}_c^* ) = 0$. From this results, we deduce that $\dim \ker { \mathcal{A}_2 } $ is exactly the index of $(\mathcal{C}_2)$ which is the index of the Dirac operator
$$\mathcal{D}^{\nabla} : \Gamma (S^{1,0} \otimes S^{5,0} ) \rightarrow \Gamma( S^{0,1} \otimes S^{5,0} ) \, . $$
By the Atiyah-Singer index theorem, 
$$\begin{array}{rcl}
index \, \mathcal{D}^{\nabla}  & = &  \{ ch ( S^{5,0} ) \hat{A} (\mathbb{S}^4) \} [\mathbb{S}^4 ]   \\
 & = & (6 + 35 ch_2(S^{1,0} ) ) (1 - p_1 /24 ) [\mathbb{S}^4]  \\
 & = & 35 \, .
\end{array}
$$
\end{proof}

\subsection{Moduli space}
In this section, we will end the proof of theorem \ref{thm2}. Here we must be more precise in our notations. If $g$ is a conformal class of metric on $\mathbb{S}^4$, there is a subbundle $S^{5,1}_g$ of $T^* \mathbb{S}^4 \otimes \Lambda ^2 T^* \mathbb{S}^4 \otimes \Lambda^2 T \mathbb{S}^4 $ associated to the representation $S^{5,1}$ and $g$. In the same way, one defines $S^{6,0}_g$ in $\Lambda^2 T^* \mathbb{S}^4 \otimes \Lambda ^2 T^* \mathbb{S}^4 \otimes \Lambda^2 T \mathbb{S}^4 $.

\begin{rem}
We have seen that each $H \in \mathcal{H}$ defines a conformal class of metrics on $\mathbb{S}^4$. In fact, the quaternionic contact structure $H$ defines a true metric on $\mathbb{S}^4$. Indeed, if we come back to section \ref{G-invariant}, the vector bundle $E$ is an oriented bundle which gives an volume form on $\Lambda^2_+ T \mathbb{S}^4$ such that $R^E$ preserves the two orientations. Then, we can choose the metric on $\mathbb{S}^4$ which gives the same volume form on $\Lambda^2_+ T \mathbb{S}^4$. We obtain a well defined map
$$ G : \mathcal{H} \rightarrow \mathcal{M} \, ,$$
where $\mathcal{M}$ is the set of smooth metrics on $\mathbb{S}^4$. The round metric on $\mathbb{S}^4$ is called $g_0$ and is the metric $G(H^{can})$.
\end{rem}

 With the help of the canonical structure $H^{can}$, we identify $\mathcal{H}$ with an open subset in $\Gamma ( T^* \mathbb{S}^4 \otimes S^{2,0}_{g_0} )$. Let $p^{i,j}$ be the orthogonal projection with respect to $g_0$ in $S^{i,j}_{g_0}$ ( $S^{i,j}_{g_0}$ will appear at most one time in our vector bundles so that the $p^{i,j}$ are well defined ). We restrict ourselves to a neighbourhood $\mathcal{U}$ of $H^{can}$ in $\mathcal{H}$ where $p^{5,1}$ (resp. $p^{6,0}$ ) gives by restriction an isomorphism from $S^{5,1}_{G(H)}$ onto $S^{5,1}_{g_0}$ ( resp. from $S^{6,0}_{G(H)}$ onto $S^{6,0}_{g_0}$ ). With the identifications given by the $p^{i,j}$, one gets maps
$$ \mathcal{T} : \mathcal{U} \rightarrow \Gamma (S^{5,1}_{g_0} ) \, , \; H \mapsto  p^{5,1}(T \! or (H))  \, , $$
 and
$$ \mathcal{B} : \mathcal{U} \oplus \Gamma(S^{5,1}_{g_0} )\rightarrow (S^{6,0}_{g_0}) \, , \; ( H , T ) \mapsto p^{6,0}(\mathcal{B}_H (T )) \, . $$

 Because of the Bianchi identity of lemma \ref{bianchi}, we have $\mathcal{B}( H, \mathcal{T} (H) ) =0$ for $ H \in \mathcal{U}$. We want to apply an implicit function theorem so we must work in Banach spaces. We assume now that our sections are $C^{k+2,\alpha}$ (H\"older-spaces). We have seen in section 6.3 that we can search a slice in $\Gamma(S^{3,1} )$. We put $\mathcal{U}_1 = \mathcal{U} \cap \Gamma(S^{3,1} ) $. Let us define the smooth map
$$ \Psi :\mathcal{U}_1 \rightarrow \mbox{Im}\, \tilde{\mathcal{D}}^* \oplus \ker \mathcal{B}_c \, , \;  a  \mapsto ( \tilde{\mathcal{D}}^* ( a) , p \, \mathcal{T}(a)  ) \, ,   $$
where $p$ is the projection on $\ker \mathcal{B}_c$ in the direction of $\mbox{Im}\, \mathcal{B}_c^* $. Because of the vanishing of the second homology group of $(\mathcal{C})$, the differential $d_{H^{can}} \psi$ is surjective. Its kernel is $\ker ( \tilde{\mathcal{D}}^* \oplus \mathcal{A} )$ and is of finite dimension $35$. Therefore, there is a submanifold $X^{35} \subset \ker \tilde{\mathcal{D}}^*\oplus \Gamma (S^{6,0})$ such that on a neighbourhood of $H^{can} $ in $\mathcal{U}_1 $, one has $\Psi(a) = 0 $ iff $a \in X^{35}$. Because of the vanishing of the homology groups $H^2$ and $H^3$, we can apply the inverse function theorem with the Bianchi operator $\mathcal{B}$ at $(H^{can}, 0)$ in order to obtain that if $p \mathcal{T}(a) =0$ then $\mathcal{T}(a) = 0$ for $a$ sufficiently small. We obtain a neighbourhood $V$ of $H^{can}$ such that
$$( \tilde{\mathcal{D}}^*(a), \mathcal{T}(a)) =0 \; \mbox{iff} \; a \in M = X^{35} \cap V \, . $$
We have obtained a $35$-dimensional family of integrable $C^{k+2, \alpha}$ quaternionic contact structures on $\mathbb{S}^7$. If $a \in M$, it satisfies a non-linear but elliptic equation, hence $a$ is smooth.

 The isotropy group $G$ of $H^{can}$ under the action of $\mathcal{G}$ is $Sp(2)$. Because $\ker \tilde{\mathcal{D}}^* $ is $Sp(2)$-invariant and $Sp(2)$ is compact, we can assume that $M$ is stable under the action of $Sp(2)$. Hence, the manifold $M$ is not the moduli space of integrable quaternionic contact structures. Nevertheless, the only diffeomorphisms acting on $M$ are in $Sp(2)$. It follows from the properness of the action of $\mathcal{G}$ on $\mathcal{H}$: an element $\phi \in \mathcal{G}$ gives a diffeomorphism $\psi$ on $\mathbb{S}^4$ acting on the metrics $G(\mathcal{H})$. The diffeomorphism $\phi$ is determined up a gauge transformation by $\psi$. The both nice behaviours of the action of diffeomorphisms on the metrics and of the gauge transformations on the connection
s give the properness of the action of $\mathcal{G}$.

Therefore there exists a neighbourhood of $[H^{can}]$ in $\mathcal{H} / \mathcal{G}$ which is homeomorphic to a neighbourhood of $H^{can}$ in $M$ quotiented by $Sp(2)$. It gives the theorem \ref{thm2}, and using the theorem 0.4 of \cite{Biq02}, one gets the corollary \ref{cor2}.

\bigskip

Among these, there is a family obtained as the boundary of quaternionic quotient constructed by Galicki in \cite{Gal91}. Let us describe these more precisely. Choose $D \in \mathfrak{sp}(2)$ and let $S^D$ be
$$S^D  = \{ x \in \mathbb{H}^k, |x|^2 + \frac{|x^* D x |^2}{4} =1 \} \, . $$
Here $x^*$ means the adjoint of $x$ with respect to the canonical quaternionic hermitian metric of $\mathbb{H}^2$. $S^D$ is isomorphic to the $7$-sphere and invariant under the diagonal action of $Sp(1)$ on right. One has the codimension $3$-distribution
$$H_x^D = \{ v \in \mathbb{H}^2, x^* v - \frac{x^* Dx }{4} ( x^* D v + v^* D x ) = 0 \} \subset T_x S^D \, . $$

This is a quaternionic contact structure which is the conformal infinity of an $AQH$ quaternionic-K\"ahler metric on the interior $B^D$ of $S^D$. Therefore $H^D$ is an integrable quaternionic contact structure. Remark that $H^D_x$ is different from the subspace of $T_x S^D \subset \mathbb{H}^2$ stable under the right-action of $\mathbb{H}$. 
The isotropy group of $H^D$ is a quotient of $K \times Sp(1)$ where $K$ is the subgroup of elements of $Sp(2)$ which commute with $D$.




\subsection{Concluding remarks} We have shown that an integrable quaternionic contact distribution on $\mathbb{S}^7$ close to the canonical one is the conformal infinity of a quaternionic-K\"ahler metric on the ball $B^8$.

 A quaternionic K\"ahler manifold can be defined with the help of a parallel $4$-form $\Omega$ with stabilizer $Sp(n)Sp(1)$. Swann \cite{Swa89} showed that in dimension greater than $8$, if $\Omega$ is closed, then $\Omega$ is parallel. On the other hand, one can construct an $8$-manifold with closed $\Omega$ which is not parallel, \cite{Sal01}. So one can ask if a quaternionic contact structure in dimension $7$  is the conformal infinity of an asymptotically hyperbolic metric associated to a closed $4$-form with stabilizer $Sp(2)Sp(1)$.


\begin{thebibliography}{Bro-Die}


\bibitem[Bes87]{Bes87}  A. L. Besse, \textit{Einstein manifolds},
Springer-Verlag, Berlin, 1987.

\bibitem[Biq00]{Biq00}  O. Biquard, \textit{M\'etriques d'Einstein
asymptotiquement sym\'etriques}, Ast\'erisque, 265 (2000).

\bibitem[Biq02]{Biq02} O. Biquard, \textit{M\'etriques autoduales sur la boule}, Invent. Math., 148, 545-607 (2002).


\bibitem[Gal91]{Gal91}  K. Galicki, \textit{\ Multi-centre metrics with negative cosmological constant }, Class. Quantum Grav., 8 (1991), 1529-1543.

\bibitem[Gra91]{Gra91} C. Robin Graham and John M. Lee, \textit{Einstein metrics with prescribed conformal infinity on the ball}, Advances in Math. 87, 186-225 (1991). 

\bibitem[Kon]{Kon}  M. Konishi, \textit{\ On manifolds with Sasakian
3-structures over quaternion-K{\"a}hlerian manifolds}, Kodai  Math. Sem.
Reps., 26 (1975), 194-200.

\bibitem[LeB82]{LeB82} C. LeBrun, \textit{$\mathcal{H}$-Space with a cosmological constant}, Proc. R. Soc. Lond. A 380, 171-185 (1982).

\bibitem[LeB89]{LeB89}  C. LeBrun, \textit{Quaternionic-K\"ahler manifolds
and conformal geometry}, Math. Ann. 284 (1989), p. 353-376.

\bibitem[LeB91]{LeB91}  C. LeBrun, \textit{On complete quaternionic-K\"ahler
manifolds}, Duke Math. J. 63 (1991), p. 723-743.

\bibitem[Mont02]{Mont02}  R. Montgomery, \textit{A tour of subriemannian
geometries, their geodesics and applications}, Math. Surv. Mon. 91 (2002),
AMS.



\bibitem[Sal82]{Sal82} S. Salamon, \textit{ Quaternionic K\"ahler manifolds}, Invent. math., 67(1982), p. 143-171.

\bibitem[Sal01]{Sal01} S. Salamon, \textit{Almost parallel structures}, Global differential geometry: the mathematical legacy of Alfred Gray (Bilbao, 2000), p.162-181.


\bibitem[Swa89]{Swa89} A. Swann, \textit{Aspects symplectiques de la g\'eom\'etrie quaternionique }, C. R. Acad. Sci. Paris, t. 308, S\'erie I, p. 225-228, 1989.

\end{thebibliography}
\end{document}